\documentclass[sts]{imsart}

\RequirePackage{amsgen,amsmath,amstext,amsbsy,amsopn,amssymb,amscd,graphicx,cancel,needspace}

\usepackage{titlesec}
\titlespacing\section{0pt}{12pt plus 2pt minus 2pt}{5pt}
\titlespacing\subsection{0pt}{10pt plus 2pt minus 2pt}{3pt}

\startlocaldefs

\newcommand{\TS} {\textstyle}

\newcommand{\Reals} {{\mathrm{I}\!\mathrm{R}}}
\newcommand{\0} {{\boldsymbol{0}}}
\newcommand{\zero} {{\boldsymbol{0}}}

\renewcommand{\P} {{\boldsymbol{P}}}
\newcommand{\Q} {{\boldsymbol{Q}}}
\newcommand{\Yvec} {{Y}}
\newcommand{\yvec} {{y}}
\newcommand{\Yveci} {{Y_i}}

\newcommand{\Phat} {{\boldsymbol{\widehat{P}}}}
\newcommand{\PhatN} {{\Phat_{\!\!N}}}
\newcommand{\PhatX} {{\Phat_{\!\!\Xvec}}}
\newcommand{\PhatYandX} {{\Phat_{\!\!\Yvec\!\!,\Xvec}}}

\newcommand{\PbootM} {{\boldsymbol{P_{\!M}^*}}}
\newcommand{\pt} {{p^{w}}}
\newcommand{\Pt} {{\boldsymbol{{P}}^{w}}}
\newcommand{\PtYandX} {{\P_{\!\!\Yvec\!\!,\Xvec}^{w}}}
\newcommand{\PtYgivenX} {{\P_{\!\!\Yvec\!|\Xvec}^{w}}}
\newcommand{\PtX} {{\P_{\!\!\Xvec}^{w}}}
\newcommand{\E} {{\boldsymbol{E}}}
\newcommand{\EP} {{\boldsymbol{E}_{\!\P}}}
\newcommand{\Ehat} {\boldsymbol{\widehat{E}}}
\newcommand{\V} {{\boldsymbol{V}}}
\newcommand{\VP} {{\boldsymbol{V}_{\!\!\!\P}}}

\newcommand{\IC} {\boldsymbol{I\!F}}

\newcommand{\AV} {{\boldsymbol{A\!V}}}

\newcommand{\X} {{\boldsymbol{X}}}

\newcommand{\Y} {{\boldsymbol{Y}}}

\newcommand{\yhat} {{\hat{y}}}

\newcommand{\cP} {{\cal P}}
\newcommand{\cR} {{\cal R}}
\newcommand{\Dist} {{\cal D}}
\newcommand{\Norm} {{\cal N}}
\newcommand{\cN} {{\cal N}}
\newcommand{\sig} {{\sigma}}

\newcommand{\rmse}{{m}}

\newcommand{\cL} {{\cal L}}

\newcommand{\Bbeta} {{\boldsymbol{\beta}}}

\newcommand{\Bbetahat} {{\boldsymbol{\hat{\beta}}}}

\newcommand{\thetahat} {{\hat{\theta}}}
\newcommand{\sigmahat} {{\hat{\sigma}}}

\newcommand{\yi} {{y_i}}

\newcommand{\argmin}{\mathrm{argmin}}

\newcommand{\xvec} {{\vec{\boldsymbol{x}}}}
\newcommand{\xvect} {{\vec{{\boldsymbol{\xi}}}}}
\newcommand{\xveci} {{\vec{\boldsymbol{x}}_i}}
\newcommand{\Xvec} {{\vec{\boldsymbol{X}}}}
\newcommand{\Xveci} {{\vec{\boldsymbol{X}}_{\!i}}}
\newcommand{\Xvecone} {{\vec{\boldsymbol{X}}_{\!\!1}}}
\newcommand{\XvecN} {{\vec{\boldsymbol{X}}_{\!\!N}}}

\newcommand{\cX} {{\cal X}}
\newcommand{\cY} {{\cal Y}}

\newcommand{\Btheta} {{\boldsymbol{\theta}}}
\newcommand{\Bthetahat} {{\boldsymbol{\hat{\theta}}}}
\newcommand{\BthetahatN} {{\boldsymbol{\hat{\theta{}}_{\!N}}}}
\newcommand{\BthetaX} {{\boldsymbol{\theta}(\X)}}
\newcommand{\BthetaP} {{\boldsymbol{\theta}(\P)}}

\newcommand{\BthetabootM} {{{\boldsymbol{\theta}_{\!M}^*}}}
\newcommand{\Bthetabootb} {{{\boldsymbol{\theta}_{\!b}^*}}}

\newcommand{\BTheta} {{\boldsymbol{\Theta}}}
\newcommand{\Bthetadot} {{\Btheta(\cdot)}}
\newcommand{\Bpsi} {{\boldsymbol{\psi}}}

\newcommand{\BLambda} {{\boldsymbol{\Lambda}}}

\newcommand{\BOmega} {{\boldsymbol{\Omega}}}

\newcommand{\xy}  {{x\textrm{-}y}}

\newcommand{\Tr}  {{'}}

\newcommand{\med} {{\mathrm{median}}}

\newcommand{\PofX} {{\P_{\!\!\Xvec}}}
\newcommand{\PofXprime} {{\P'_{\!\!\Xvec}}}
\newcommand{\PofYgivenX} {{\P_{\!\!\Yvec\!|\Xvec}}}
\newcommand{\PofYgivenXi} {{\P_{\!\!\Yvec\!|\Xveci}}}

\newcommand{\PofYgivenx} {{\P_{\!\!\Yvec\!|\Xvec=\xvec}}}
\newcommand{\PofYandX} {{\P_{\!\!\Yvec,\Xvec}}}

\newcommand{\PofuYgivenX} {{\P_{\!\!Y\!|\Xvec}}}

\newcommand{\Qmodel} {{\Q_{\Yvec\!|\Xvec;\Btheta}}}
\newcommand{\Qmodelx} {{\Q_{\Yvec\!|\Xvec=\xvec;\Btheta}}}
\newcommand{\Qmodelzero} {{\Q_{\Yvec\!|\Xvec;\Btheta_0}}}
\newcommand{\Qmodelxzero} {{\Q_{\Yvec\!|\Xvec=\xvec;\Btheta_0}}}

\newcommand{\KL}{{D_{K\!L}}}

\newcommand{\tw}{{\tilde{w}}}
\newcommand{\hw}{{\hat{w}}}

\endlocaldefs

\begin{document}


\begin{frontmatter}

  \title{ Models as Approximations II: A Model-Free Theory of
    Parametric Regression} \runtitle{Models as Approximations II}

\begin{aug}

  \author{\fnms{Andreas}
          \snm{Buja}\thanksref{t1,m1}\ead[label=e1]{buja.at.wharton@gmail.com}},
  \author{\fnms{Lawrence}
          \snm{Brown}\thanksref{t1,m1}\ead[label=e3]{lbrown@wharton.upenn.edu}},
  \author{\fnms{Arun}
          \snm{Kumar Kuchibhotla}\thanksref{m1}\ead[label=e3]{arun@wharton.upenn.edu}},
  \author{\fnms{Richard}
          \snm{Berk}\thanksref{m1}\ead[label=e2]{berkr@wharton.upenn.edu}},
  \author{\fnms{Ed}
          \snm{George}\thanksref{t3,m1}\ead[label=e3]{egeorge@wharton.upenn.edu}},
  \and
  \author{\fnms{Linda}
          \snm{Zhao}\thanksref{t1,m1}\ead[label=e3]{lzhao@wharton.upenn.edu}},

  \thankstext{t1}{Supported in part by NSF Grant DMS-10-07657 and  and DMS-1310795.} 
  \thankstext{t3}{Supported in part by NSF Grant DMS-14-06563.}

  \runauthor{A. Buja et al.}

  \affiliation{The Wharton School -- University of Pennsylvania\thanksmark{m1} }

  \address{Statistics Department,
           The Wharton School, University of Pennsylvania,
           400 Jon M. Huntsman Hall, 3730 Walnut Street,
           Philadelphia, PA 19104-6340
           \printead{e1}.}

\end{aug}

\maketitle

\begin{abstract}
  We develop a model-free theory of general types of parametric
  regression for iid observations.  The theory replaces the parameters
  of parametric models with statistical functionals, to be called
  ``regression functionals'', defined on large non-parametric classes
  of joint $\xy$ distributions, without assuming a correct model.
  Parametric models are reduced to heuristics to suggest plausible
  objective functions.  An example of a regression functional is the
  vector of slopes of linear equations fitted by OLS to largely
  arbitrary $\xy$ distributions, without assuming a linear model (see
  Part~I).  More generally, regression functionals can be defined by
  minimizing objective functions or solving estimating equations at
  joint $\xy$ distributions.  In this framework it is possible to
  achieve the following: (1)~define a notion of well-specification for
  regression functionals that replaces the notion of correct
  specification of models, (2)~propose a well-specification diagnostic
  for regression functionals based on reweighting distributions and
  data, (3)~decompose sampling variability of regression functionals
  into two sources, one due to the conditional response distribution
  and another due to the regressor distribution interacting with
  misspecification, both of order $N^{-1/2}$, (4)~exhibit
  plug-in/sandwich estimators of standard error as limit cases of
  $\xy$ bootstrap estimators, and (5)~provide theoretical heuristics
  to indicate that $\xy$ bootstrap standard errors may generally be
  more stable than sandwich estimators.
\end{abstract}

\begin{keyword}[class=AMS]
\kwd[Primary ]{62J05}      
\kwd{62J20}                
\kwd{62F40}                
\kwd[; secondary ]{62F35}  
\kwd{62A10}                
\end{keyword}

\begin{keyword}
\kwd{Ancillarity of regressors}
\kwd{Misspecification}
\kwd{Econometrics}
\kwd{Sandwich estimator}
\kwd{Bootstrap}
\kwd{Bagging}
\end{keyword}

\end{frontmatter}


\newpage

{\em
``The hallmark of good science is that it uses models and 'theory' but never believes them.''
(J.W.~Tukey, 1962\nocite{ref:Tukey-1962}, citing Martin Wilk)
}

\bigskip

\section{Introduction}
\label{sec:intro}

We develop in this second article a model-free theory of parametric
regression, assuming for simplicity iid $\xy$ observations with quite
arbitrary joint distributions.  The starting point is the realization
that regression models are approximations and should not be thought of
as generative truths.  A general recognition of this fact may be
implied by the commonly used term ``working model,'' but this vague
term does not resolve substantive issues, created here by the fact
that models are approximations and not truths.  The primary issue is
that model parameters define meaningful quantities only under
conditions of model correctness.  If the idea of models as
approximations is taken seriously, one has to extend the notion of
parameter from model distributions to basically arbitrary
distributions.  This is achieved by what is often called ``projection
onto the model,'' that is, finding for the actual data distribution
the best approximating distribution within the model; one defines that
distribution's parameter settings to be the target of estimation.
Through such ``projection'' the parameters of a working model are
extended to ``statistical functionals,'' that is, mappings of largely
arbitrary data distributions to numeric quantities.  We have thus
arrived at a functional point of view of regression, a view based on
what we call {\bf\em regression functionals}.

The move from traditional regression parameters in correctly specified
models to regression functionals obtained from best approximations may
raise fears of opening the gates to irresponsible data analysis where
misspecification is of no concern.  No such thing is intended here.
Instead, we rethink the essence of regression and develop a new notion
of {\bf\em well-specification of regression functionals}, to replace
the notion of {\em correct specification of regression models}.  In
the following bullets we outline an argument in the form of simple
postulates.
\begin{itemize}
\item The essence of regression is the asymmetric analysis of
  association: Variables with a joint distribution $\P$ are divided
  into response and regressors.
\item Motivated by prediction and causation problems, interest focuses
  on properties of the conditional distribution of the response given
  the regressors.
\item The goal or, rather, the hope is that the chosen
  quantities/functionals of interest are properties of the observed
  conditional response distribution, irrespective of the regressor
  distribution.
\item Consequently, a regression functional will be called {\bf\em
    well-specified} if it is a property of the observed conditional
  response distribution at hand, {\em irrespective of the regressor
    distribution}.
\end{itemize}
The first bullet is uncontroversial: asymmetric analysis is often
natural, as in the contexts of prediction and causation.  The second
bullet remains at an intended level of vagueness as it explains the
nature of the asymmetry, namely, the focus on the
regressor-conditional response distribution.  Intentionally there is
no mention of regression models.  The third bullet also steers clear
of regression models by addressing instead quantities of interest,
that is, regression functionals.  In this and the last bullet, the
operational requirement is that the quantities of interest not depend
on the regressor distribution.  It is this constancy across regressor
distributions that turns the quantities of interest into
properties of the conditional response distribution alone.

All this can be made concrete with reference to the groundwork laid in
Part~I, Section 4.  Consider the regression functional consisting of
the coefficient vector obtained from OLS linear regression.  It was
shown in Part~I that this vector does not depend on the regressor
distribution (is well-specified) if and only if the conditional
response mean is a linear function of the regressors.
Thus the coefficient vector fully describes the conditional mean
function, but no other aspect of the conditional response
distribution.  Well-specification of the OLS coefficient functional is
therefore a weaker condition than correct specification of the linear
model by setting aside homoskedasticity and Gaussianity which are
linear model requirements not intimately tied to the slopes.

A desirable feature of the proposed definition of well-specification
is that it generalizes to arbitrary types of parametric regression or,
more precisely, to the statistical functionals derived from them.  In
particular, it applies to GLMs where the meaning of well-specified
coefficients is again correct specification of the mean function but
setting aside other model requirements.  Well-specification further
applies to regression functionals derived from optimizing general
objective functions or solving estimating equations.
Well-specification finally applies to any ad hoc quantities if they
define regression functionals for joint $\xy$ distributions.

The proposed notion of well-specification of regression functionals
does not just define an ideal condition for populations but also lends
itself to a tangible methodology for real data.  A diagnostic for
well-specification can be based on perturbation of the regressor
distribution without affecting the conditional response distribution.
Such perturbations can be constructed by reweighting the joint $\xy$
distribution with weight functions that only depend on the regressors.
If a regression functional is not constant under such reweighting, it
is misspecified.

In practice, use of this diagnostic often works out as follows.  Some
form of misspecification will be detected for some of the quantities
of interest, but the diagnostic will also aid in interpreting the
specifics of the misspecification.  The reason is that reweighting
essentially localizes the regression functionals.  For the
coefficients of OLS linear regression, for example, this means that
reweighting reveals how the coefficients of the best fitting linear
equation vary as the weight function moves across regressor space.
Put this way, the diagnostic seems related to non-parametric
regression, but its advantage is that it focuses on the quantities of
interest at all times, while switching from parametric to
non-parametric regression requires a rethinking of the meaning of the
original quantities in terms of the non-parametric fit.  To guide
users of the diagnostic to insightful choices of weight functions, we
introduce a set of specific reweighting methodologies, complete with
basic statistical inference.

Following these methodological proposals, we return to the inferential
issues raised in Part~I and treat them in generality for all types of
well-behaved regression functionals.  We show that sampling variation
of regression functionals has two sources, one due to the conditional
response distribution, the other due to the regressor distribution
interacting with misspecification, where ``misspecification'' is meant
in the sense of ``violated well-specification'' of the regression
functional.  A central limit theorem (CLT) shows that {\em both}
sources, as a function of the sample size $N$, are of the usual order
$N^{-1/2}$.  Finally, it is shown that asymptotic plug-in/sandwich
estimators of standard error are limits of $\xy$ bootstrap estimators,
revealing the former to be an extreme case of the latter.

The present analysis becomes necessarily more opaque because algebra
that worked out explicitly and lucidly for linear OLS in Part~I is
available in the general case only in the form of asymptotic
approximation based on influence functions.  Still, the analysis is
now informed by the notion of well-specification of regression
functionals, which gives the results a rather satisfactory form.

The article continues as follows.  In Section~\ref{sec:targets} we
discuss typical ways of defining regression functionals, including
optimization of objective functions and estimating equations.  In
Section~\ref{sec:well-specification} we give the precise definition of
well-specification and illustrate it with various examples.  In
Section~\ref{sec:reweighting} we introduce the reweighting diagnostic
for well-specification, illustrated in
Section~\ref{sec:reweighting-methodology} with specific reweighting
methodologies applied to the LA homeless data (Part~I).
Section~\ref{sec:decomposition} shows for plug-in estimators of
regression functionals how the sampling variability is canonically
decomposed into contributions from the conditional response noise and
from the randomness of the regressors.  In Section~\ref{sec:CLT} we
state general CLTs analogous to the OLS versions of Part~I.  In
Section~\ref{sec:plugin-bootstrap} we analyze model-free estimators of
standard error derived from the $M$-of-$N$ pairs bootstrap and
asymptotic variance plug-in (often of the sandwich form).  It holds in
great generality that plug-in is the limiting case of bootstrap when
$M \rightarrow \infty$.  We also give some heuristics to suggest that
boostrap estimators might generally be more stable than
plug-in/sandwich estimators.  In Section~\ref{sec:summary} we
summarize the path taken in these two articles.

{\bf Remark}: For notes on the history of model robustness, see
Part~I, Section~1.  For the distinction between model robustness and
outlier/heavy-tail robustness, see Part~I, Section~13.


\section{Targets of Estimation: Regression Functionals}
\label{sec:targets}

This section describes some of the ways of constructing regression
functionals, including those based on ``working models'' used as
heuristics to suggest plausible objective functions.  We use the
following notations and assumptions throughout: At the population
level there are two random variables, the regressor $\Xvec$ with
values in a measurable space $\cX$ and the response $\Yvec$ with
values in a measurable space $\cY$, with a joint distribution
$\PofYandX$, a conditional response distribution $\PofYgivenX$ and a
marginal regressor distribution~$\PofX$.  We express the connection
between them using ``$\otimes$'' notation:
\begin{equation} \label{eq:notation}
\PofYandX ~=~ \PofYgivenX \otimes \PofX .
\end{equation}
Informally this is expressed in terms of densities by
$p(\yvec,\xvec) = p(\yvec|\xvec) p(\xvec)$.  In contrast to Part~I,
the regressor and response spaces $\cX$ and $\cY$ are now entirely
arbitrary.  The typographic distinction between $\Xvec$ and $\Yvec$ is
a hold-over from the OLS context of Part~I.  Both spaces, $\cX$ {\em
  and} $\cY$, can be of any measurement type, univariate or
multivariate, or even spaces of signals or images.

Regression functionals need to be defined on universes of joint
distributions that are sufficiently rich to grant the manipulations
that follow, including the assumed existence of moments, influence
functions, and closedness for certain mixtures.  The details are
tedious, hence deferred to Appendix~\ref{sec:assumptions} without
claim to technical completeness.  The treatment is largely informal so
as not to get bogged down in distracting detail.  Also, the
asymptotics will be traditional in the sense that $\cX$ and $\cY$ are
fixed and $N \rightarrow \infty$.  For more modern technical work on
related matters, see Kuchibhotla et al.~(2018)\nocite{ref:Kuchi-2018}.


\subsection{Regression Functionals from Optimization: ML and PS Functionals}
\label{sec:ML-PS-functionals}

In Part~I we described the interpretation of linear OLS coefficients
as regression functionals.  The expression ``linear OLS'' is used on
purpose to avoid the expression ``linear models'' because no model is
assumed.  Fitting a linear equation using OLS is a procedure to
achieve a best fit of an equation by the OLS criterion.  This approach
can be generalized to other objective functions
$\cL(\Btheta;\yvec,\xvec)$:
\begin{equation} \label{eq:objective-fct}
  \Btheta(\P) ~=~ \argmin_{\,\Btheta \in \BTheta} \, \EP[\cL(\Btheta;\Yvec,\Xvec)]
\end{equation}
A common choice for $\cL(\Btheta;y,\xvec)$ is the negative
log-likelihood of a parametric regression model for $\Yvec|\Xvec$,
defined by a parametrized family of conditional response distributions
$\{ \Qmodel\!\!: \Btheta \!\in\! \BTheta \}$ with conditional
densities
$\{ q(\yvec\,|\,\xvec;\Btheta)\!\!: \Btheta \!\in\! \BTheta \}$.  The
model is not assumed to be correctly specified, and its only purpose
is to serve as a heuristic to suggest an objective function:
\begin{equation} \label{eq:neg-log-lik}
\cL(\Btheta;\yvec,\xvec) ~=~ - \log \, q(\yvec\,|\,\xvec;\Btheta) .
\end{equation}
In this case the regression functional resulting from
\eqref{eq:objective-fct} will be called a maximum-likelihood
functional or {\em ML functional} for short.
It minimizes the
Kullback-Leibler (KL) divergence of
$\PofYandX = \PofYgivenX \otimes \PofX$ and $\Qmodel \otimes \PofX$,
which is why one loosely interprets an ML functional as arising from a
``projection of the actual data distribution onto the parametric
model.''  ML functionals can be derived from major classes of
regression models, including GLMs.  Technically, they also comprise
many M-estimators based on Huber $\rho$ functions (Huber
1964\nocite{ref:Huber-1964}), including least absolute deviation (LAD,
$L_1$) as an objective function for conditional medians, and tilted
$L_1$ versions for arbitrary conditional quantiles, all of which can
be interpreted as negative log-likelihoods of certain distributions,
even if these may not usually be viable models for actual data.  Not
in the class of negative log-likelihoods are objective functions for
M-estimators with redescending influence functions such as Tukey's
biweight estimator (which also poses complications due to
non-convexity).

Natural extensions of ML functionals can be based on so-called
``proper scoring rules'' (Appendix~\ref{ref:proper-scoring}) which
arise as cross-entropy terms of Bregman divergences A special case is
the expected negative log-likelihood arising as the cross-entropy term
of KL divergence.  The optimization criterion is the proper scoring
rule applied to the conditional response distribution $\PofYgivenX$
and model distributions $\Qmodel$, averaged over regressor space
with~$\PofX$.  The resulting regression functionals may be called
``proper scoring functionals'' or simply {\em PS functionals}, a
superset of ML functionals.  All PS functionals, including ML
functionals, have the important property of Fisher consistency: If the
model is correctly specified, i.e., if $\,\exists\, \Btheta_0$ such
that $\PofYgivenX = \Qmodelzero$, then the population minimizer
is~$\Btheta_0$:
\begin{equation} \label{eq:Fisher-consistency}
  \mbox{if}~~\PofYandX = \Qmodelzero \otimes \PofX, ~~\mbox{then}~~ \Btheta(\P) = \Btheta_0 .
\end{equation}
See Appendix~\ref{ref:proper-scoring} for background on proper scoring
rules, Bregman divergences, and some of their robustness properties to
outliers and heavy tailed distributions.

Further objective functions are obtained by adding parameter penalties
to existing objective functions:
\begin{equation} \label{eq:penalty}
\tilde{\cL}(\Btheta;\yvec,\xvec) ~=~ \cL(\Btheta;\yvec,\xvec) + \lambda \cR(\Btheta) .
\end{equation}
Special cases are ridge and lasso penalties.  Note that
\eqref{eq:penalty} results in one-parameter families of penalized
functionals $\Btheta_\lambda(\P)$ defined for populations as well,
whereas in practice $\lambda\!=\!\lambda_N$ applies to finite $N$ with
$\lambda_N \!\rightarrow\! 0$ as $N\rightarrow\infty$.


\subsection{Regression Functionals from Estimating Equations: EE Funtionals}
\label{sec:EE}

Objective functions are often minimized by solving stationarity
conditions that amount to estimating equations with the scores
$\Bpsi(\Btheta;\yvec,\xvec) = -\nabla_\Btheta \cL(\Btheta;\yvec,\xvec)$:
\begin{equation} \label{eq:EE}
  \EP[ \Bpsi(\Btheta;\Yvec,\Xvec) ] ~=~ \0.
\end{equation}
One may generalize and define regression functionals as solutions in
cases where $\Bpsi(\Btheta;\yvec,\xvec)$ is not the gradient of an
objective function; in particular it need not be the score function of
a negative log-likelihood.  Functionals in this class will be called
{\em EE functionals}.  For OLS, the estimating equations are
the normal equations, as the score function for the slopes is
\begin{equation} \label{eq:OLS-moments}
  \Bpsi_{OLS}(\Bbeta;y,\xvec) ~=~ \xvec y  \!-\! \xvec \xvec\Tr \, \Bbeta
  ~=~ \xvec (y  \!-\! \xvec\Tr \, \Bbeta)  .
\end{equation}
A seminal work that inaugurated asymptotic theory for general
estimating equations is by Huber~(1967)\nocite{ref:Huber-1967}.  A
more modern and rigorous treatment is in
Rieder~(1994)\nocite{ref:Rieder-1994}.

An extension is the ``Generalized Method of Moments'' (GMM, Hansen
1982\nocite{ref:Hansen-1982}).  It applies when the number of moment
conditions (the dimension of $\Bpsi$) is larger than the dimension of
$\Btheta$.  An important application is to causal inference based on
numerous instrumental variables.

Another extension is based on ``Generalized Estimating Equations''
(GEE, Liang and Zeger 1986\nocite{ref:LZ-1986}).  It applies to
clustered data that have intra-cluster dependence, allowing
misspecification of the variance and intra-cluster dependence.


\subsection{The Point of View of Regression Functionals and its Implications}
\label{sec:regression-functional-point-of-view}

Theories of parametric models deal with the issue that a traditional
model parameter has many possible estimators, as in the normal model
$\cN(\mu,\sigma^2)$ where the sample mean is in various ways the
optimal estimate of $\mu$ whereas the median is a less efficient
estimate of the same $\mu$.  The comparison of estimates of the same
traditional parameter has been proposed as a basis of misspecification
tests (Hausman~1978)\nocite{ref:Hausman-1978} and called ``test for
parameter estimator inconsistency''
(White~1982)\nocite{ref:White-1982}.  In a framework based on
regression functionals the situation presents itself differently.
Empirical means and medians, for example, are not estimators of the
same parameter; instead, they represent different statistical
functionals.  Similarly, slopes obtained by linear OLS and linear LAD
are different regression functionals.  Comparing them by forming
differences creates new regression functionals that may be useful as
diagnostic quantities, but in a model-robust framework there is no
concept of ``parameter inconsistency'' (White~1982,
p.$\,$15\nocite{ref:White-1982}), only a concept of differences between
regression functionals.  

A further point is that in a model-robust theory of observational (as
opposed to causal) association, there is no concept of ``omitted
variables bias.''  There are only regressions with more or fewer
regressor variables, none of which being ``true'' but some being more
useful or insightful than others.  Slopes in a larger regression are
distinct from the slopes in a smaller regression.  It is a source of
conceptual confusion to write the slope of the $j$'th regressor as
$\beta_j$, irrespective of what the other regressors are.  In more
careful notation one indexes slopes with the set of selected
regressors $M$ as well, $\beta_{j \cdot M}$, as is done of necessity
in work on post-selection inference (e.g., Berk et
al.~2013)\nocite{ref:Berk-et-al-2013}.  Thus the linear slopes
$\beta_{j \cdot M}$ and $\beta_{j \cdot M'}$ for the $j$'th regressor,
when it is contained in both of two regressor sets $M \neq M'$, should
be considered as distinct regression functionals.  The difference
$\beta_{j \cdot M'} - \beta_{j \cdot M}$ is not a bias but a
difference between two regression functionals.  If it is zero, it
indicates that the difference in adjustment between $M$ and $M'$ is
immaterial for the $j$'th regressor.  If $\beta_{j \cdot M'}$ and
$\beta_{j \cdot M}$ are very different with opposite signs, there
exists a case of Simpson's paradox for this regressor.

It should be noted that regression functionals generally depend on the
full joint distribution~$\PofYandX$ of the response {\em and} the
regressors.  Conventional regression parameters describe the
conditional response distribution only under correct specification,
$\PofYgivenX = \Qmodel$, while the regressor distribution $\PofX$ is
sidelined as ancillary.  That the ancillarity argument for the
regressors is not valid under misspecification was documented in
Part~I, Section~4.  In the following sections this fact will be the
basis of the notion of well-specification of regression functionals.

Finally, we state the following to avoid misunderstandings: In the
present work, the objective is not to recommend particular regression
functionals, but to point out the freedoms we have in choosing them
and the conceptual clarifications we need when using them.


\section{Mis-/Well-Specification of Regression Functionals}
\label{sec:well-specification}

The introduction motivated a notion of well-specification for
regression functionals, and this section provides the technical
notations.  The heuristic idea is that a regression functional is
well-specified for a joint distribution of the regressors and the
response if it does not depend on the marginal regressor distribution.
In concrete terms, this means that the functional does not depend on
where the regressors happen to fall.  The functional is therefore a
property of the conditional response distribution alone.


\subsection{Definition of Well-Specification for Regression Functionals}
\label{sec:define-well-specification}

Recall the notation introduced in \eqref{eq:notation}:
$\PofYandX = \PofYgivenX \otimes \PofX$.  Here a technical detail
requires clarification: conditional distributions are defined only
almost surely with regard to $\PofX$, but we will assume that
$\xvec \mapsto \P_{\!\!\Yvec|\Xvec = \xvec}$ is a Markov kernel
defined for all $\xvec \in \cX$.\footnote{Thus we assume a ``regular
  version'' has been chosen, as is always possible on Polish spaces.}
With these conventions, $\PofYgivenX$ and $\PofX$ uniquely determine
$\PofYandX = \PofYgivenX \otimes \PofX$ by \eqref{eq:notation}, but
not quite vice versa.  Thus $\Btheta(\cdot)$ can
be written as
\begin{equation*} 
  \Btheta(\P) ~=~ \Btheta(\PofYgivenX \otimes \PofX).
\end{equation*}

\medskip

\noindent{\bf Definition:~} {\em The regression functional
  $\Bthetadot$ is well-specified for $\PofYgivenX$ if
  \begin{equation*}
  \Btheta(\PofYgivenX \otimes \PofX) ~=~ \Btheta(\PofYgivenX \otimes \PofXprime)
  \end{equation*}
  for all acceptable regressor distributions $\PofX$
  and~$\PofXprime$.  }

\medskip

\noindent
The term ``acceptable'' accounts for exclusions of regressor
distributions such as those due to non-identifiability when fitting
equations, in particular, perfect collinearity when fitting linear
equations (see Appendix~\ref{sec:assumptions}).

\smallskip

\noindent {\bf Remarks:}
\begin{itemize}
\item Importantly, the notion of well-specification is a {\em joint
    property} of a specific $\Bthetadot$ and a specific $\PofYgivenX$.
  A regression functional will be well-specified for some conditional
  response distributions but not for others.
\item The notion of well-specification represents an idealization, not
  a reality.  Well-specification is never a fact, only degrees of
  misspecification are.  Yet, idealizations are useful because they
  give precision and focus to an idea.  Here, the idea is that a
  regression functional is intended to be a property of the
  conditional response distribution $\PofYgivenX$ alone, regardless of
  the regressor distribution~$\PofX$.
\end{itemize}


\subsection{Well-Specification --- Some Exercises and Special Cases}

\noindent
Before stating general propositions, here are some special cases to train
intuitions.

\begin{itemize} \itemsep 0.5em

\item The OLS slope functional can be written
  $\Bbeta(\P) = \EP[\Xvec \Xvec\Tr]^{-1} \EP[\Xvec \mu(\Xvec)]$, where
  $\mu(\xvec) = \EP[Y|\Xvec\!=\!\xvec]$.  Thus $\Bbeta(\P)$ depends
  on $\PofYgivenX$ only through the conditional mean function.  The
  functional is well-specified if $\mu(\xvec) = \Bbeta_0\Tr \xvec$ is
  linear, in which case $\Bbeta(\P) = \Bbeta_0$.  For the reverse, see
  Part~I, Proposition~4.1.

\item A special case is regression through the origin, which we
  generalize slightly as follows.  Let $h(\xvec)$ and $g(\yvec)$ be two
  non-vanishing real-valued square-integrable functions of the
  regressors and the response, respectively.  Define
  \[
  \Btheta_{h,g}(\P) ~=~ \frac{\EP[ \,g(\Yvec) h(\Xvec)\, ]}{\EP[ \,h(\Xvec)^2 ]} .
  \]
  Then $\Btheta_{h,g}(\P)$ is well-defined for $\PofYgivenX$ if
  $\EP[\,g(\Yvec) | \Xvec\,] = c \cdot h(\Xvec)$ for some~$c$.

\item An ad hoc estimate of a simple linear regression slope is
  \[
    \Btheta(\P) ~=~ \EP[ (Y' \!-\! Y'') / (X' \!-\! X'') \;\big|\big.\;|X' \!-\! X''| > \delta ],
  \]
  where $(Y',X'), (Y'',X'') \sim \P$ iid and $\delta > 0$.  It is
  inspired by Part~I, Section~10 and Gelman and
  Park~(2008)\nocite{ref:GP-2008}.  It is well-specified if
  $\EP[Y|X] = \beta_0 + \beta_1 X$, in which case
  $\Btheta(\P) = \beta_1$.

\item Ridge regression also defines a slope functional.  Let $\BOmega$
  be a symmetric non-negative definite matrix and
  $\Bbeta\Tr \BOmega \Bbeta$ its quadratic penalty.  Solving the
  penalized LS problem yields
  $\Bbeta(\P) = (\EP[\Xvec \Xvec\Tr] + \BOmega)^{-1} \EP[\Xvec
  \mu(\Xvec)]$.  This functional is well-specified if the conditional
  mean is linear, $\mu(\xvec) = \Bbeta_0\Tr \xvec$ for some
  $\Bbeta_0$, and $\BOmega = c \EP[\Xvec \Xvec\Tr]$ for some
  $c \!\ge\! 0$, in which case $\Bbeta(\P) = 1/(1\!+\!c)\, \Bbeta_0$,
  causing uniform shrinkage across all regression coefficients.

\item Given a univariate response $Y$, what does it mean for the
  functional $\Btheta(\P) \!=\! \EP[\,Y]$ to be well-specified
  for $\PofuYgivenX$?  It looks as if it did not depend on the
  regressor distribution and is therefore always well-specified.  This
  is a fallacy, however.  Because $\EP[Y] \!=\! \EP[\mu(\Xvec)]$,
  it follows that $\EP[Y]$ is independent of $\PofX$ iff the
  conditional response mean is constant: $~\mu(\Xvec) = \EP[Y]$.

\item Homoskedasticity: The average conditional variance functional
  $\sig^2(\P) = \EP[\VP[Y|\Xvec]]$ is well-specified iff
  $\VP[Y|\Xvec=\xvec] \!=\! \sig_0^2$ is constant, in which case
  $\sig^2(\P) = \sig_0^2$.  A difficulty is that access to this
  functional assumes a correctly specified mean function
  $\mu(\Xvec) = \EP[Y|\Xvec]$.  

\item The average conditional MSE functional wrt linear OLS is
  $\E[(Y\!-\Bbeta(\P)'\Xvec)^2] = \E[\rmse^2(\Xvec)]$ using the
  notation of Part~I.  If it is well-specified, that is, if
  $\rmse^2(\Xvec) = \rmse_o^2$ is constant, then linear model-based
  inference is asymptotically justified (Part~I, Lemma~11.4~(a)).

\item The correlation coefficient $\rho(Y,X)$, if interpreted as a
  regression functional in a regression of $Y$ on $X$, is
  well-specified only in the trivial case when $\mu(X)$ is constant
  and $\VP[Y] > 0$, hence $\rho(Y,X) = 0$.

\item Fitting a linear equation by minimizing least absolute
  deviations (LAD, the $L_1$ objective function) defines a regression
  functional that is well-specified if there exists $\Bbeta_0$ such
  that $~\med[\,\PofuYgivenX\,] = \Bbeta_0\Tr \Xvec$.

\item In a GLM regression with a univariate response and canonical
  link, the slope functional is given by
  \[
    \Bbeta(\P) ~=~
    \argmin_{\Bbeta} \, \EP \big[\, b\,(\Xvec\Tr\Bbeta) - Y \Xvec\Tr\Bbeta \,\big] ,
  \]
  where $b(\theta)$ is a strictly convex function on the real line and
  $\theta = \xvec\Tr \Bbeta$ is the ``canonical parameter'' modeled by
  a linear function of the regressors.  The stationary equations
  are\footnote{To avoid confusion with matrix transposition, we write
    $\partial b$ instead of $b'$ for derivatives.}
  \[
    \EP \big[ Y \Xvec \big] ~=~
    \EP \big[\, \partial b\,(\Xvec\Tr\Bbeta) \Xvec \,\big] .
  \]
  This functional is well-specified iff
  $\EP \big[ Y | \Xvec \big] = \partial b \,(\Xvec\Tr\Bbeta)$
  for~$\Bbeta=\Bbeta(\P)$.  Well-specification of $\Bbeta(\P)$ has
  generally no implication for $\VP \big[ Y | \Xvec \big]$, except
  in the next example.

\item Linear logistic regression functionals are a special case of GLM
  functionals where $Y \!\in\! \{0,1\}$ and
  $b(\theta) = \log(1+\exp(\theta))$.  Well-specification holds iff
  $\P[\,Y\!=\!1\,|\Xvec] = \phi(\Xvec'\Bbeta)$ for~$\Bbeta=\Bbeta(\P)$
  and $\phi(\theta) = \exp(\theta)/(1+\exp(\theta))$.  Because the
  conditional response distribution is Bernoulli, the conditional mean
  of $Y$ determines the conditional response distribution uniquely,
  hence well-specification of the regression functional $\Bbeta(\P)$
  is the same as correct specification of the logistic regression
  model.

\item If $\Btheta(\P)$ is well-specified for $\PofYgivenX$, then so is
  the functional $f(\Btheta(\P))$ for any function $f(\cdot)$.  An
  example in linear regression is the predicted value
  $\Bbeta(\P)\Tr \xvec$~ at the regressor location~$\xvec$.  Other
  examples are contrasts such as $\beta_1(\P) \!-\! \beta_2(\P)$ where
  $\beta_j(\P)$ denotes the $j$'th coordinate of $\Bbeta(\P)$.
\item A meaningless case of ``misspecified functionals'' arises when
  they do not depend on the conditional response distribution at all:
  $\Btheta(\PofYgivenX \otimes \PofX) = \Btheta(\PofX)$.  Examples
  would be tabulations and summaries of individual regressor
  variables.  They could not be well-specified for~$\PofYgivenX$
  unless they are constants.

\end{itemize}


\subsection{Well-Specification of ML, PS and EE Functionals}
\label{sec:well-spec-ML-PS-EE}

The following lemma, whose proof is obvious, applies to all ML
functionals.  The principle of pointwise optimization in regressor
space covers also all PS functionals (see
Appendix~\ref{sec:PS-for-regression},
equation~\eqref{eq:pointwise-min}).

\medskip

\noindent{\bf Proposition \ref{sec:well-spec-ML-PS-EE}.1:} {\em If
  $\Btheta_0$ minimizes
  $\EP[\cL(\Yvec|\Xvec;\Btheta)\,|\,\Xvec\!=\!\xvec]$ for all
  $\xvec \!\in\! \cX$, then the minimizer $\Btheta(\P)$ of
  $\EP[\cL(\Yvec|\Xvec;\Btheta)]$ is well-specified for $\PofYgivenX$,
  and $\Btheta(\PofYgivenX \otimes \PofX) = \Btheta_0$ for all
  acceptable regressor distributions~$\PofX$.  }

\medskip

\noindent
The following fact is corollary of
Proposition~\ref{sec:well-spec-ML-PS-EE}.1 but could have been gleaned
from Fisher consistency \eqref{eq:Fisher-consistency}.

\medskip

\noindent{\bf Proposition \ref{sec:well-spec-ML-PS-EE}.2:} {\em If
  $\Bthetadot$ is a ML or PS functional for the working model
  $\{ \Qmodel\!: \Btheta\!\in\!\BTheta \}$, it is well-specified for
  all model distributions $\PofYgivenX = \Qmodel$.  }

\medskip

\noindent
The next fact states that an EE functional is well-specified for
a conditional response distribution if it satisfies the EE
conditionally and globally across regressor space for one
value~$\Btheta_0$.

\medskip

\noindent{\bf Proposition \ref{sec:well-spec-ML-PS-EE}.3:} {\em If
  $\Btheta_0$ solves
  $\EP[ \Bpsi(\Btheta_0;\Yvec,\Xvec) | \Xvec\!=\!\xvec ] = \0$ for all
  $\xvec \in \cX$, then the EE functional defined by
  $\EP[ \Bpsi(\Btheta;\Yvec,\Xvec) ] = \0$ is well-specified for
  $\PofYgivenX$, and $\Btheta(\PofYgivenX \otimes \PofX) = \Btheta_0$
  for all acceptable regressor distributions~$\PofX$.  }

\medskip

\noindent The proof is in Appendix~\ref{sec:proof}.


\subsection{Well-Specification and Causality}
\label{sec:causality}

The notion of well-specification for regression functionals relates to
aspects of causal inference based on direct acyclic graphs (DAGs) and
the Markovian structures they represent (e.g., Pearl
2009)\nocite{ref:Pearl-2009}.  Given a DAG, the theory explains which
choices of regressors $\Xvec$ permit correct descriptions of causal
effects for a given outcome variable $\Yvec$.  Focusing on one such
choice of $\Xvec$ and $\Yvec$, one is left with the task of describing
interesting quantitative aspects of the conditional distribution
$\PofYgivenX$, which is thought to be unchanging under different
manipulations and/or sampling schemes of the regressors $\Xvec$.
Therefore, if a quantity of interest is to describe causal effects
properly, it should do so irrespective of where the values of the
causal variables $\Xvec$ have fallen.  This is exactly the requirement
of well-specification for regression functionals.  In summary, proper
causal effects must arise as quantities of interest that are
well-specified in the sense of
Section~\ref{sec:define-well-specification}.

Recently, Peters, B\"uhlmann and Meinshausen (2016,
Section~1.1)\nocite{ref:PBM-2016} discussed a related notion of
``invariance'' which can be interpreted as ``invariance to regressor
distributions.''  They propose this notion as a heuristic for causal
discovery and inference based on multiple data sources with the same
variables, one variable being singled out as the response $\Yvec$.
These multiple data sources are leveraged as follows: If for a subset
of variables, $\Xvec$, the association $\Xvec \rightarrow \Yvec$ is
causal, then the conditional distribution $\PofYgivenX$ will be the
same across data sources.  Subsets of causal variables $\Xvec$ with
shared $\PofYgivenX$ across sources may therefore be discoverable if
the sources differ in their regressor distributions and/or
interventions on causal variables.  For concreteness, the authors
focus on a linear structural equation model (SEM), which allows us to
reinterpret their proposals by abandoning the SEM assumption and
consider instead the regression functional consisting of the OLS
regression coefficients resulting from the linear SEM.  Thus the
proposed method is at heart an approach to detecting and inferring
well-specified quantities, cast in a causal framework.

In the following section we will introduce a diagnostic for
well-specification that can be interpreted as emulating multiple data
sources from a single data source.  The proposal is to systematically
reweight the data to synthetically create alternative datasets.
Peters et al.~(2016, Section~3.3)\nocite{ref:PBM-2016} briefly mention
the idea of conditioning as related to the idea of multiple data
sources.  Such conditioning is naturally achieved by locally
reweighting the data, as will be shown next.


\medskip

\section{A Reweighting Diagnostic for Well-Specification: Targets
  and Inferential Tools}
\label{sec:reweighting}

Well-specification of regression functionals connects naturally to
reweighting, both of populations and of data.  A concrete illustration
of the basic idea can be given by again drawing on the example of
linear OLS: The OLS slope functional is well-specified iff
$\EP[\Yvec|\Xvec] = \Bbeta_0' \Xvec$ for some $\Bbeta_0$, in
which case for any non-negative weight function $w(\xvec)$ we have
$\Bbeta_0 = \argmin_\Bbeta \, \EP[w(\Xvec) \, (\Yvec - \Bbeta'
\Xvec)^2]$.  Therefore the reweighting of interest is with regard to
weights that are functions of the regressors only.  The general reason
is that such weights affect the distribution of the regressors but not
the conditional response distribution.  Reweighting provides an
intuitive basis for diagnosing well-specification of regression
functionals.  Because of the practical importance of the proposed
reweighting diagnostic, we insert this material early, deferring
estimation and inference to Section~\ref{sec:decomposition}.


Reweighting has an extensive history in statistics, too rich to
recount.  The present purpose of reweighting is methodological: to
diagnose the degree to which the null hypothesis of well-specification
of a regression functional is violated.  To this end we propose what
we call a ``tilt test.''  It provides evidence of whether a
real-valued regression functional is likely to rise or fall (tilt up
or down) from one extreme of reweighting to another.  The conclusions
from a rejection based on this test are simple and interpretable.

In practice, the majority of regression functionals of interest are
regression slopes connected to specific regressors.  A more
interesting problem than detection of misspecification is another
question: Does misspecification impinge on the statistical
significance of a slope of interest?  That is, would a slope have lost
or gained statistical significance if the regressor distribution had
been different?  This is the primary question to be addressed by the
reweighting diagnostic.


\subsection{Reweighting and Well-Specification}
\label{sec:reweighting-and-well-specification}

Consider reweighted versions of the joint distribution
$\P = \PofYandX$ with weight functions $w(\xvec)$ that depend only on
the regressors, not the response, written as
\[
  \PtYandX(d\yvec,d \xvec) \!=\! w(\xvec) \PofYandX(d\yvec,d\xvec) ,
  ~~~\mbox{or}~~~
  \pt(\yvec,\xvec) = w(\xvec) \, p(\yvec,\xvec) ,
\]
where $w(\xvec) \!>\! 0$ and $\EP[w(\Xvec)] \!=\! 1$, which turns
$\PtYandX$ into a joint probability distribution for $(\Yvec,\Xvec)$
with the same support as $\PofYandX$.  At times, for specific weight
functions, we will write $w(\Xvec) \PofYandX$ instead of $\PtYandX$.

\medskip

\noindent {\bf Lemma \ref{sec:reweighting-and-well-specification}}:~~$\PtYgivenX = \PofYgivenX$ ~and~ $\PtX = w(\Xvec) \PofX$.

\medskip

\noindent
The proof is elementary and simplest in terms of densities:
\begin{eqnarray*}
&& \pt(\xvec) = {\TS \int} \pt(\yvec,\xvec) d\yvec = {\TS \int} w(\xvec)
p(\yvec,\xvec) d\yvec = w(\xvec) {\TS \int} p(\yvec,\xvec) d\yvec = w(\xvec) p(\xvec) ,
\\
&& \pt(\yvec|\xvec) = \pt(\yvec,\xvec)/ \pt(\xvec) =
(w(\xvec)\,p(\yvec,\xvec)) / (w(\xvec)\,p(\xvec) = p(\yvec,\xvec)/p(\xvec) =
p(\yvec|\xvec).
\end{eqnarray*}
We obtain as an immediate consequence:

\medskip

\noindent {\bf Proposition
  \ref{sec:reweighting-and-well-specification}}:~{\em If the
  regression functional $\Bthetadot$ is well-specified for
  $\PofYgivenX$, it is unchanged under arbitrary $\Xvec$-dependent
  reweighting:~ $\Btheta(\PtYandX) = \Btheta(\PofYandX)$.}

\medskip

\noindent {\bf Remark:} In fixed-$\X$ linear models theory, which
assumes correct specification, it is known that reweighting the data
with fixed weights grants unbiased estimation of coefficients.
Translated to the current framework, this fact returns as a statement
of invariance of well-specified functionals under $\Xvec$-dependent
reweighting.

\smallskip

Tests of misspecification based on reweighting were proposed by White
(1980a, Section~4)\nocite{ref:White-1980a} for linear OLS.  The
approach generalizes to arbitrary types of regression and regression
functionals as follows: Given a weight function $w(\Xvec)$ normalized
for $\P$, the null hypothesis is ~$H_0: \Btheta(\Pt) =\Btheta(\P)$.
For the case that $\Btheta(\cdot)$ is the vector of OLS linear
regression coefficients, White (ibid.,
Theorem~4)\nocite{ref:White-1980a} proposes a test statistic based on
plug-in estimates and shows its asymptotic null distribution to be
$\chi^2$.  The result is a Hausman test
(1978)\nocite{ref:Hausman-1978} whereby (using model-oriented
language) an efficient estimate under the model is compared to an
inefficient but consistent estimate.  Rejection indicates
misspecification.  We will not draw on White's results but instead use
the $\xy$ bootstrap as a basis of inference because (1)~it directly
applies to general types of regression under mild technical
conditions, and (2)~it lends itself to augmentation of visual displays
that provide more informative diagnostics than vanilla tests.  White
(1980a)\nocite{ref:White-1980a} did not develop a methodology for
reweighting tests other than recommending experimentation with
multiple weight functions.  The present goal is to introduce highly
interpretable {\em one-parameter families of weight functions} and to
illustrate their practical use to gain insights into the nature of
misspecifications.


\subsection{The Well-Specification Diagnostic: Population Version}
\label{sec:reweighting-diagnostic}

In order to construct interpretable weight functions, we construct
them as functions of a {\em univariate variable} $Z$.  This variable
will often be one of the real-valued regressors, $Z\!=\!X_j$.
However, the variable $Z$ may be any function of the regressors,
$Z \!=\! f(\Xvec)$, as when $Z \!=\! \Bbeta'\Xvec$ is the OLS fit of
$Y$, or $Z \!=\! X_{j\bullet}$ is $X_j$ adjusted for all other
regressors (Part~I, Section~9).\footnote{Mathematically, the
  restriction to weights as a function of univariate variables $Z$ is
  no restriction at all because any $w(\xvec)$ can be trivially
  described as the identity function of $Z = w(\cdot)$. }

Given a variable $Z$, consider for concreteness a univariate Gaussian
weight function of $Z$, centered at~$\xi$ on the $Z$~axis:
\begin{equation} \label{eq:Gaussian-weights}
w_{\xi}(z) \,=\, w_{\xi}^*(z) / \E[w_{\xi}^*(Z)] \,,
~~~~
w_{\xi}^*(z) \,\propto\, \exp(-(z-\xi)^2/(2 \gamma^2)) \,,
\end{equation}
where $\gamma$ is a user-specified bandwidth parameter (see
Section~\ref{sec:reweighting-diagnostic-data-version} below).

Next consider a one-dimensional regression functional $\theta(\P)$,
such as a linear regression slope.  A graphical diagnostic is obtained
by plotting $\theta(\cdot)$ as a function of the reweighting
centers~$\xi$:
\begin{equation} \label{eq:reweighting-trace}
\xi ~\mapsto~ \theta_{\xi}(\P) ~=~ \theta(w_{\xi}(Z) \, \P) \,.
\end{equation}
If the regression functional $\theta(\P)$ is well-specified for
$\PofYgivenX$, then $\theta_{\xi}(\P)$ is constant in~$\xi$ and equal
to $\theta(\P)$.  Equivalently, if $\theta_{\xi}(\P)$ is not constant
in $\xi$, then $\theta(\P)$ is misspecified.  Thus non-constancy is a
sufficient criterion for misspecification.  Insightful choices of
traces of the form \eqref{eq:reweighting-trace} will be proposed
below.


\subsection{The Reweighting Diagnostic: Data Version}
\label{sec:reweighting-diagnostic-data-version}

To make the diagnostic actionable on data, one obtains estimates
\[
\thetahat_{\xi} ~=~ \theta(\hw_\xi(Z) \Phat),
\]
where $\hw_\xi(x)$ is a weight function that is empirically normalized
to unit mass, $\Ehat[\hw_\xi(Z)] \!=\! 1$, where $\Ehat[...]$ denotes
the sample average.  This means using weights for the observations of
the form
\[
w_i = \hw_\xi(z_i) \propto \exp(-(z_i-\xi)^2/(2 \gamma^2)), ~~~
\TS \frac{1}{N} \sum_i w_i = 1, ~~~i\!=\!1,...,N.
\]
We parametrize the bandwidth $\gamma = \alpha \sigmahat(Z)$ in terms
of the empirical standard deviation $\sigmahat(Z)$ of $Z$ and a
multiplier $\alpha$.  In the examples we use $\alpha\!=\!1$.

In order to plot a discretized version of the trace
$\xi \mapsto \thetahat_{\xi}$, we obtain estimates $\thetahat_{\xi}$
for a grid of values $\xi_{(1)}<...<\xi_{(K)}$ on the $Z$ axis, a
simple choice being the interior deciles of the empirical $Z$
distribution.  Hence $K\!=\!9$, unless $Z$ has numerous ties, causing
some deciles to collapse.  Finally, we plot
$\xi_{(k)} \mapsto \thetahat_{\xi_{(k)}}$.  This is carried out in
Figures~\ref{fig:reweighting-LA-coef-eq-PercVacant}-\ref{fig:LA-Reweighting-weights-eq-yhat}
for the LA homeless data (see
Section~\ref{sec:reweighting-methodology}).


\subsection{Interpretations of the Reweighting Diagnostic}
\label{sec:reweighting-diagnostics-interpretation}

The reweighting diagnostic is likely to be accessible to practitioners
of regression.  One reason is that the restriction to weights as a
function of a univariate variable $Z$ permits a simple left-to-right
comparison: Is $\xi \mapsto \theta(w_{\xi}(Z) \,\P)$ higher or lower
on the right than on the left?  In our experience, the dominant
feature of such traces is indeed monotonicity.  The intuitive appeal
of reweighting is further helped by two mutually compatible
interpretations:
\begin{itemize}
\item {\bf Data frequency}: Reweighting mimics scenarios of datasets
  that contain more or fewer observations as a function of~$Z$ than
  the observed dataset.  Thus it answers questions such as ``what if
  there were more observations with low (or high) values of~$Z$?''  In
  this sense reweighting mimics alternative data sources based on the
  data at hand.
\item {\bf Conditioning}: Reweighting can be seen as ``soft
  conditioning on $Z$'' in the sense that conditioning on ``sharp
  inclusion'' in an interval $\xi \!-\!c \!<\! Z \!<\! \xi \!+\! c$ is
  replaced by ``soft inclusion'' according to the weight
  function~$w_{\xi}(z)$.  In this sense reweighting localizes the
  regression functional.  However, note that when $Z = X_j$, for
  example, the localization is of ``codimension 1'' in regressor space
  (approached as the bandwidth $\gamma \rightarrow 0$).
\end{itemize}
In what follows we use either of these interpretations depending on
the context.


\subsection{Inferential Features for Reweighting Diagnostics}

\begin{figure}[!th]
  \centering
  \includegraphics[width=5.5in]{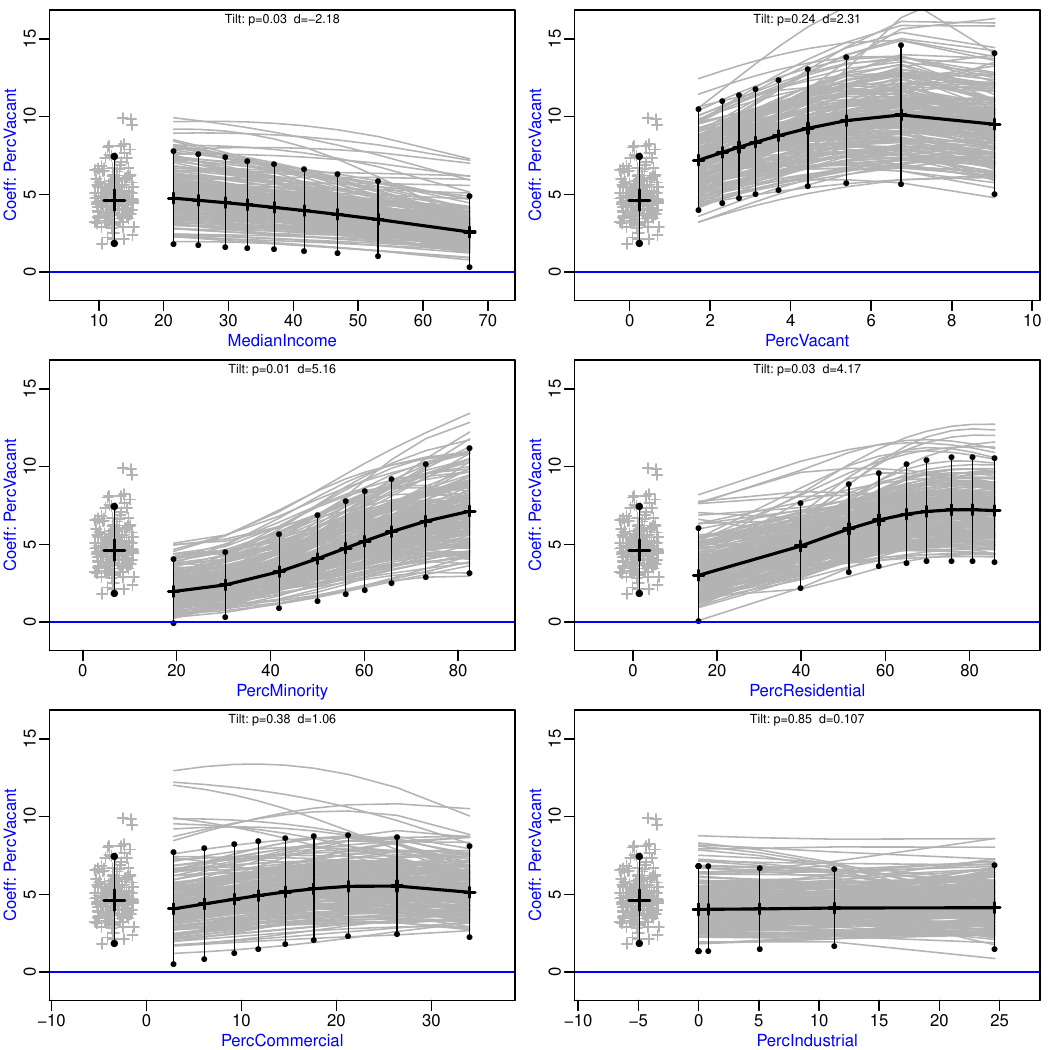}
  \vspace{-.1in}
  \caption{\it Diagnostics for the slope of {\tt PercVacant}; LA
    Homeless Data (see Section 2, Part~I). \newline Vertical axis =
    regression coefficient of {\tt PercVacant} (in all frames);
    horizontal axes = regressors.  If the vertical axis is interpreted
    causally as the effect size of {\tt PercVacant} on the response
    {\tt StreetTotal} (of homeless in a census tract), the following
    can be inferred: the effect size of {\tt PercVacant} is greater
    for high values of {\tt PercMinority} (center left frame) and low
    values of {\tt MedianInc} (top left frame), and possibly also for
    high values of {\tt PercResidential}.  Near the top margin of each
    frame are the p-values of a ``Tilt'' test for the difference
    between the right-most and left-most effect sizes.}
\nocite{ref:BKY-2008}
  \label{fig:reweighting-LA-coef-eq-PercVacant}
\end{figure}

Graphical diagnostics need inferential augmentation to answer
questions of whether visually detected features are real.  Presently
the two main questions are:
\begin{itemize} \itemsep 0em
\item[(1)] Is the variation/non-constancy in
  $\xi_{(k)} \mapsto \thetahat_{\xi_{(k)}}$ sufficiently strong to be
  statistically significant and hence suggest misspecification
  of~$\theta(\cdot)$?
\item[(2)] Where are the estimates $\thetahat_{\xi_{(k)}}$
  statistically significantly different from zero?
\end{itemize}
For regression slopes question (2) may be more relevant than (1)
because one usually cares about their statistical significance.
Therefore, to answer question (2), we decorate the diagnostic plot
with traces of bootstrapped estimates, as shown in the plots of
Figures~\ref{fig:reweighting-LA-coef-eq-PercVacant}-\ref{fig:LA-Reweighting-weights-eq-yhat}.
Bootstrap resampling is done from the actual, not the reweighted,
data.  The weight functions have the same centers $\xi_{(k)}$, but
their bandwidth is based on bootstrapped standard deviations.  In the
figures we show 199 bootstrap traces in gray color, amounting to a
so-called ``spaghetti plot''.  Along with the bootstrap replications
we also show bootstrap error bars at the grid locations.  Their widths
are $\pm 2$ bootstrap standard errors.

As can be illustrated with
Figures~\ref{fig:reweighting-LA-coef-eq-PercVacant}-\ref{fig:LA-Reweighting-weights-eq-yhat},
statistical significance can feature 
a variety of patterns.  Significance may exist~...
\begin{itemize}
\item[(2a)] ...~across the whole range of reweighting centers
  $\xi_{(k)}$ and in the same direction, as in the top right plot of
  Figure~\ref{fig:reweighting-LA-coef-eq-PercVacant};
\item[(2b)] ...~both on the left and the right but in opposite
  directions with a transition through insignificance in between, as
  is nearly the case in the center left plot of
  Figure~\ref{fig:LA-Reweighting-weights-eq-regressors};
\item[(2c)] ...~over part of the range, typically the left or the
  right side; such tendencies are seen in the two center plots of
  Figure~\ref{fig:reweighting-LA-coef-eq-PercVacant};
\item[(2d)] ...~nowhere, as in the bottom right plot of
  Figure~\ref{fig:LA-Reweighting-weights-eq-regressors}.
\end{itemize}

To answer question~(1) regarding the presence of misspecification, we
piggyback on the bootstrap exercise meant to answer question~(2).
Because most detections of misspecification arise from a monotone tilt
in the trace $\xi_{(k)} \mapsto \thetahat_{\xi_{(k)}}$, we construct a
cheap test statistic by forming the difference between the two extreme
points of the trace,
$\thetahat_{\xi_{(K)}} - \thetahat_{\xi_{(1)}}$.\footnote{This test
  statistic does not result in a Hausman
  (1978)\nocite{ref:Hausman-1978} test: both estimates are
  ``inefficient under correct model specification.''  However, it
  quantifies an obvious visual feature of the traces.}  We obtain its
bootstrap distribution almost for free, hence we can perform a crude
bootstrap test by placing the null value zero in the bootstrap
distribution.  The bootstrap p-value and the test statistic are shown
near the top of each plot frame in
Figures~\ref{fig:reweighting-LA-coef-eq-PercVacant}-\ref{fig:LA-Reweighting-weights-eq-yhat}.
For example, the top left frame of
Figure~\ref{fig:reweighting-LA-coef-eq-PercVacant} shows ``{\tt
  Tilt:$\;$p=0.04~d=2.18}'', meaning that the difference of $2.18$ is
statistically significant with a (two-sided) p-value
of~0.04.\footnote{For 199 bootstrap replicates the lowest possible two-sided
  p-value is~0.01 $=2\cdot1/(1+199)$.}

Finally we show on the left side of each frame a visual version of
unweighted plain statistical significance of the quantity of interest
in the form of a bootstrap confidence interval around the unweighted
estimate $\thetahat \pm 2$ unweighted bootstrap standard errors.  In
addition, we show 199 bootstrap estimates (gray points horizontally
jittered to reduce overplotting).  The location on the horizontal axis
has no meaning other than being 10\% to the left of the range
$(\xi_{(1)},\xi_{(K)})$ of the traces.


\section{The Reweighting Diagnostic for Well-Specification:
  Methodology and Examples}
\label{sec:reweighting-methodology}

The following subsections demonstrate three different purposes of the
diagnostic.  The quantities of interest are linear OLS slopes, though
the approach generalizes to all types of regression that permit
reweighting:
\begin{itemize}
\item {\bf Focal Slope}: Expose a slope $\beta_k(\P)$ of special
  interest to reweighting on each regressor in turn: $Z \!=\! X_j$ for
  $j\!=\!1,...,p$ (Section~\ref{sec:focal-coefficient}).  This
  produces highly interpretable insights into interactions of
  regressor $X_k$ with all other regressors $X_j$, without modeling
  these interactions directly.
\item {\bf Nonlinearity detection}: Expose each regression slope
  $\beta_j(\P)$ to reweighting on its own regressor, $Z \!=\!X_j$
  (Section~\ref{sec:reweighting-on-own-regressor}).  This produces
  insights into marginal nonlinear behaviors of response surfaces.
\item {\bf Focal Reweighting Variable}: Use a single reweighting
  variable of interest (here: $Z = \Bbeta'\Xvec$) to diagnose
  well-specification of all components of a regression functional,
  here: slopes $\beta_j(\P)$ (Section~\ref{sec:focal-reweighting}).
\end{itemize}
These diagnostics will be illustrated with the LA homeless data of
Part~I, Section~2.  The observations consist of a sample of 505 census
tracts in the LA metropolitan area, and the variables are seven
quantitative measures of the tracts with largely self-explanatory
names: The response is the {\tt StreetTotal} (count) of homeless
people in a census tract, and the six regressors are: {\tt
  MedianIncome} (of households, in \$1,000s), {\tt PercMinority}, and
the prevalences of four types of lots: {\tt PercCommercial}, {\tt
  PercVacant}, {\tt PercResidential} and {\tt PercIndustrial}.


\subsection{Diagnostics for a Focal Regression Coefficient of Interest (Figure~\ref{fig:reweighting-LA-coef-eq-PercVacant})}
\label{sec:focal-coefficient}

One variable stands out as potentially accessible to intervention by
public policies: {\tt PercVacant}.  Vacant lots could be turned into
playgrounds, sports fields, parks, or offered as neighborhood
gardens.\footnote{Such programs have indeed been enacted in some
  cities.  We abstain from commenting on the controversies surrounding
  such policies.}  It would therefore be of interest to check whether
the regression coefficient of {\tt PercVacant} possibly measures a
causal effect, for which it is a necessary condition that it be
well-specified (Section~\ref{sec:causality}).  To this end,
Figure~\ref{fig:reweighting-LA-coef-eq-PercVacant} shows diagnostics
for the coefficient of {\tt PercVacant} under reweighting on all six
regressors.

As the plots show, statistical significance of the coefficient of {\tt
  PercVacant} holds by and large under reweighting across the ranges
of all six regressors.  While this is comforting, there exists a
weakening of significance in the extremes of the ranges of three
regressors: high {\tt MedianInc}, low {\tt PercMinority} and low {\tt
  PercResidential}.  With these qualitative observations it is already
indicated that well-specification of the coefficient of {\tt
  PercVacant} is doubtful, and indeed the tilt tests show statistical
significance with 2-sided p-values of 0.01 and 0.02 for {\tt
  PercMinority} and {\tt MedianInc}, respectively.  The variable {\tt
  PercResidential} also looks rather steep, but its tilt test has a
weaker p-value around~0.1.  Finally, a very weak indication is shown for
larger effects at higher levels of {\tt PercVacant}.

Does this indication of misspecification invalidate a causal effect of
{\tt PercVacant}?  It does not.  It only points to the likely
possibility that the causal effect is not correctly described by a
single linear regression coefficient; it is rather a more complex
function of the regressors.  Useful insight into the nature of the
causal effect (if this is what it is) can be gleaned from the
diagnostic plots by using them to answer an obvious question: Where is
the effect of {\tt PercVacant} likely to be strong?  An answer might
indeed help in prioritizing interventions.  Interpreting the plots of
Figure~\ref{fig:reweighting-LA-coef-eq-PercVacant} liberally, one
could state that the effect of {\tt PercVacant} looks strongest for
census tracts with high {\tt PercMinority}, followed by high {\tt
  PercResidential} and low {\tt MedianInc}.  These observations seem
rather plausible and may indeed point to census tracts worth
prioritizing for intervention with public policies.\footnote{An
  application of this type of diagnostic to the Boston Housing data is
  in Appendix~\ref{sec:reweighting-Boston}.}

The insights gained so far point to the presence of interactions
between {\tt PercVacant} and other regressors because the slope of
{\tt PercVacant} varies at different levels of those other regressors.
A natural next step would be more detailed modeling that includes
interactions between {\tt PercVacant} and the three interacting
regressors, but the essential insights have already been gained.


\begin{figure}[!t]
  \centering
  \includegraphics[width=5.5in]{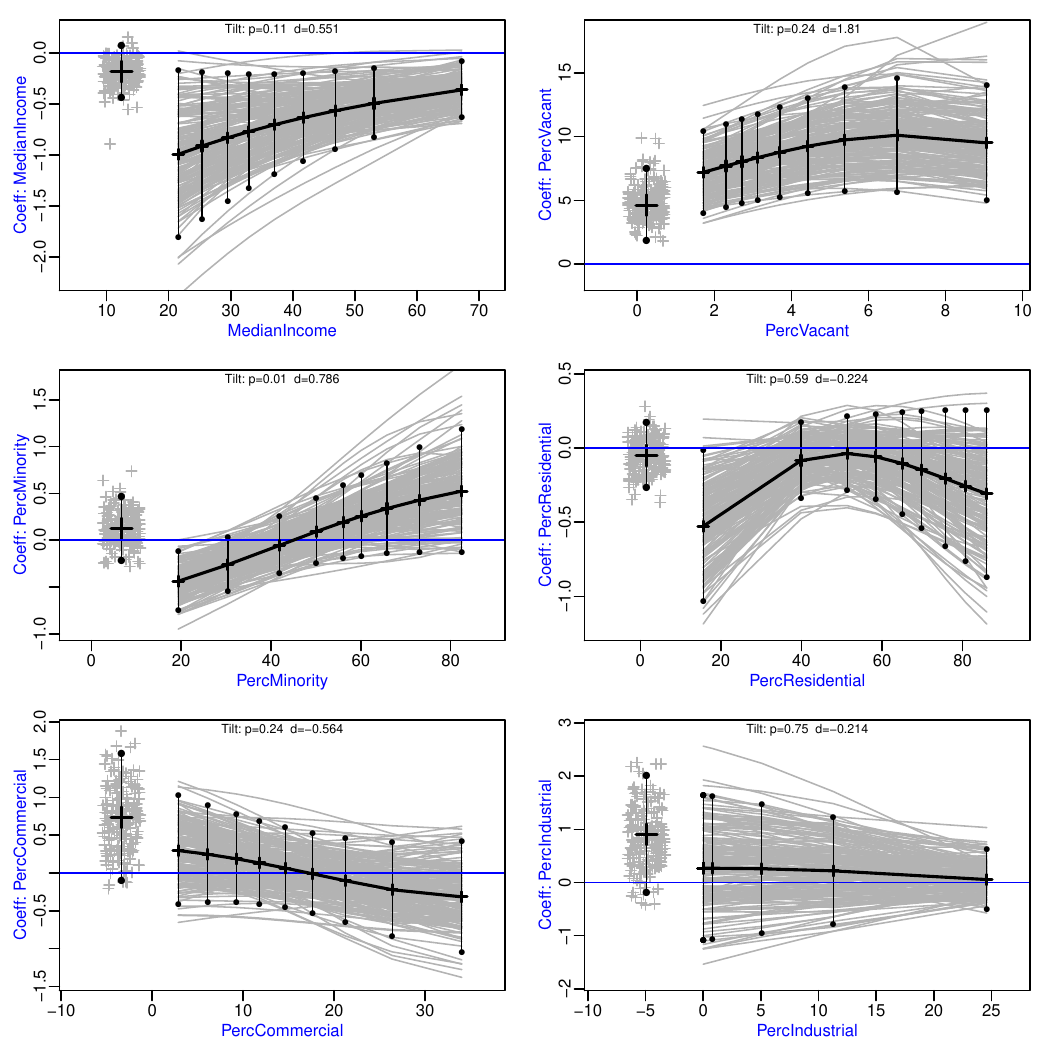}
  \vspace{-.1in}
  \caption{\it Misspecification Diagnostics: Slopes reweighted by
    their own regressors --- indications of nonlinearity.  The center
    left frame suggests that the slope of {\tt PercMinority} reverses
    from a negative slope for low values of {\tt PercMinority} to a
    positive slope for high values of {\tt PercMinority}.
  }\nocite{ref:BKY-2008}
  \label{fig:LA-Reweighting-weights-eq-regressors}
\end{figure}

\subsection{Diagnostics for Slopes Reweighted by Their Own Regressors~(Figure~\ref{fig:LA-Reweighting-weights-eq-regressors})}
\label{sec:reweighting-on-own-regressor}

The top right plot in
Figure~\ref{fig:reweighting-LA-coef-eq-PercVacant} is a special case
where the slope of interest is reweighted by its own regressor, {\tt
  PercVacant}.  It has a different interpretation, not related to
interactions but to nonlinear effects.  To get a better picture of the
possibilities that can arise in real data, we show in
Figure~\ref{fig:LA-Reweighting-weights-eq-regressors} the
corresponding plots for all six regressors and their slopes.

Glancing at the six plots, we note some unpredictable effects of
reweighting, both on the values and the estimation uncertainties of
the slopes.  We find examples of larger and smaller estimates as well
as stronger and weaker statistical significances relative to their
unweighted analogs:
\begin{itemize} \itemsep 0em
\item Bottom left plot for the regressor {\tt PercCommercial}: The
  unweighted estimate of $\beta_j(\P)$ (on the left side of the plot)
  is weakly statistically significant (the lower end of the $\pm2$
  standard error confidence interval touches zero).  The reweighted
  estimates of $\beta_j(w_{\xi}(X_j)\P)$, however, are closer to zero
  and nowhere statistically significant for any $\xi$ in the range of
  {\tt PercCommercial}.
\item Top right plot for the regressor {\tt PercVacant}: The
  unweighted estimate and the reweighted estimates are all
  statistically significant, but the reweighted ones are
  systematically larger and much more statistically significant.
\end{itemize}

Another noteworthy case of a different nature appears for the
regressor {\tt PercMinority}
(Figure~\ref{fig:LA-Reweighting-weights-eq-regressors}, center left
plot).  While the unweighted estimate is statistically insignificant,
the locally reweighted estimates reveal a striking pattern:
\begin{itemize} \itemsep 0em
\item For low values of {\tt PercMinority} $\approx 20\%$, the slope
  is negative and statistically significant: Incrementally more
  minorities is associated with a lower {\tt StreetTotal} of homeless.
\item For high values of {\tt PercMinority} $\approx 80\%$, the slope
  is positive and (weakly) statistically significant: Incrementally
  more minorities is associated with a higher {\tt StreetTotal} of
  homeless.
\end{itemize}
This finding (if real) represents a version of Simpson's paradox: In
aggregate, there is no statistically significant association, but,
conditional on low and high values of {\tt PercMinority}, there is,
and in opposite directions.

In Appendix~\ref{sec:additive-models} we discuss some reasons for the
unpredictable behaviors of slopes under reweighting wrt to their own
regressors.  We also mention a (weak) link to partial additive models
with one nonlinear term.


\begin{figure}[!t]
  \centering
  \includegraphics[width=5.5in]{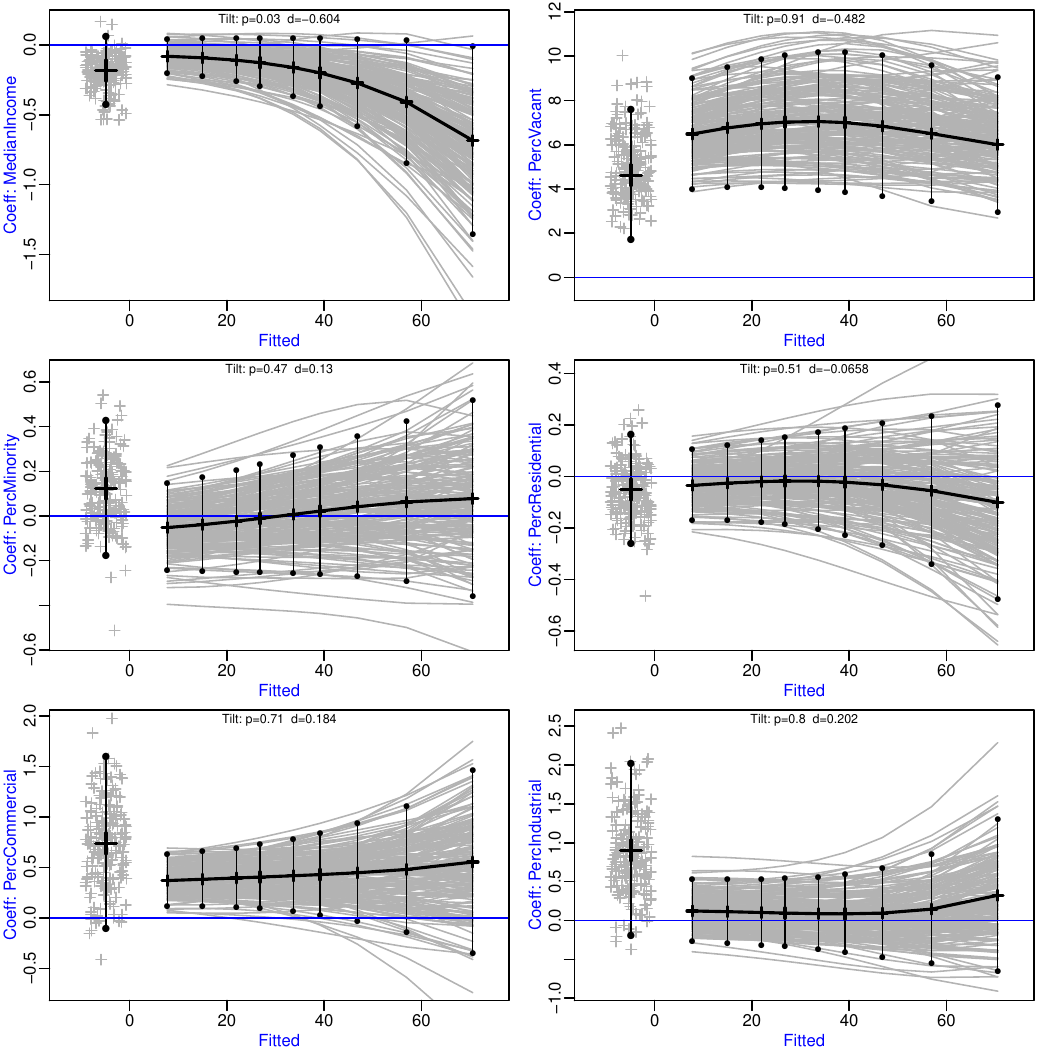}
  \vspace{-.1in}
  \caption{\it Misspecification Diagnostics using one focal
    reweighting variable, the best linear approximation/prediction
    {\tt Fitted}, for all slopes.  The fan shapes left to right
    suggest that all slope estimates except the one for {\tt
      PercVacant} have more sampling variability for higher fitted
    values.}
  \nocite{ref:BKY-2008}
  \label{fig:LA-Reweighting-weights-eq-yhat}
\end{figure}

\subsection{Diagnostics for a Focal Reweighting Variable of Interest (Figure~\ref{fig:LA-Reweighting-weights-eq-yhat})}
\label{sec:focal-reweighting}

Next we illustrate a version of the diagnostics that subjects all
slopes of a linear regression to a single reweighting variable of
interest.  The goal is to detect misspecification in any coefficient,
and the hope is to do so by reweighting based on a variable~$Z$ that
is both powerful and interpretable.  Taking a cue from traditional
residual diagnostics, we choose the OLS best approximation,
$Z = \Bbeta\Tr \Xvec$.  The data version is based on reweighting as
function of the OLS estimates of the fitted values,
$z_i = \yhat_i = \Bbetahat\Tr \xvec_i$.\footnote{The estimated slope
  vector $\Bbetahat$ is frozen across bootstraps, ignoring a
  lower-order source of sampling variability.}  The question is
whether any coefficient reveals misspecification when comparing it on
data with more low versus more high values of the linear
approximation.  The expectation is that the gradient of the linear
approximation should be a direction of high average response variation
and hence may have a higher promise of revealing misspecifications
than other directions in regressor space.

Figure~\ref{fig:LA-Reweighting-weights-eq-yhat} shows this diagnostic
applied to the LA homeless data, labeling the reweighting variable as
{\tt Fitted}.  Some observations are as follows:
\begin{itemize}
\item The only slope with signs of misspecification is for {\tt
    MedianIncome} (top left plot), whose tilt test has a p-value of
  0.03.  This slope achieves mild statistical significance for high
  values of {\tt Fitted}, which would indicate that the ``effect'' (if
  any) of differences in {\tt MedianIncome} matter more for high
  values of the linear prediction {\tt Fitted}.
\item The slope of {\tt PercCommercial} (bottom left plot) shows no
  signs of misspecification, but it is mildly statistically
  significant only for low values of {\tt Fitted} due to the
  lower estimation uncertainty in that range.
\item Five of the six plots feature a fan shape of the bootstrap
  spaghetti bands (exception: {\tt PercVacant}).  This indicates that
  these five slope estimates have greater estimation uncertainty for
  higher values of {\tt Fitted}.
\end{itemize}
The last point illustrates that the diagnostic is not only informative
about the average level of estimates but also about their estimation
uncertainty.


\subsection{Summary Comments on Reweighting Diagnostics}

The reweighting diagnostics proposed here are not meant to replace
other types of diagnostics, typically based on residual analysis.
They are, however, able to answer questions about quantities of
interest and effects of regressors that residual analysis might not.
They may also be able to provide insights into the nature of
nonlinearities and interactions without explicitly modeling them.
Furthermore they are easily augmented with inferential features such
as bootstrap spaghetti bands and tests of misspecification with
specific interpretations.  Finally, they are able to localize regions
in regressor space with high or low estimation uncertainty.


\section{Estimation of Regression Functionals:   ~~~~~~~~~~ Canonical Decomposition of Estimation Offsets}
\label{sec:decomposition}

We return to the task of building a general framework of plug-in
estimation of regression functionals based on iid data.  We decompose
sampling variability into its two sources, one due to the conditional
response distribution, the other due to the randomness of the
regressors interacting with misspecification.  Along the way we find
new characterizations of well-specification of regression functionals.


\subsection{Regression Data and Plug-In Estimation}

We adopt some of the notations and assumptions from Part~I, Section~5:
Data consist of $N$ iid draws
$(\Yveci,\Xveci) \!\sim\! \P \!=\! \PofYandX$; the responses $\Yveci$
are collected in a data structure $\Y \!=\! \{\Yveci\}_i$, and the
regressors $\Xveci$ in another data structure $\X \!=\! \{\Xveci\}_i$,
called ``data frame'' in programming languages such as R
(2008)\nocite{ref:R-2008}.  We avoid the terms ``vector'' and
``matrix'' because in a general theory of regression all variables ---
responses and regressors --- can be of any type and of any
dimension.\footnote{Recall that the typographic difference between
  $\Yvec$ and $\Xvec$ is a holdover from Part~I, where the response
  was assumed univariate quantitative.}  This is why not only $\X$ but
$\Y$ is best thought of as a (random) ``data frame.''  Regression of
$\Y$ on $\X$ is any attempt at estimating aspects of the conditional
distribution~$\PofYgivenX$.
We limit ourselves to regression functionals $\Bthetadot$ that allow
plug-in estimation $\Bthetahat = \Btheta(\Phat)$ where
$\Phat = \PhatYandX = (1/N)\sum \delta_{(\Yveci,\Xveci)}$
is the joint empirical distribution.  If necessary we may write
$\PhatN$ for~$\Phat$ and $\BthetahatN$ for~$\Bthetahat$.  In addition,
we will also need the empirical regressor distribution
$\PhatX = (1/N)\sum \delta_{\Xveci}$.


\subsection{The Conditional Parameter of Model-Trusting Fixed-$\X$ Regression}
\label{sec:conditional-parameter}

We now define the important notion of a ``conditional parameter'' for
arbitrary regression functionals, thereby providing the target of
estimation for fixed-$\X$ theories.  For OLS slopes this target of
estimation is $\Bbeta(\X) \!=\! \EP[\Bbetahat|\X]$ (Part~I,
Section~5).  We use the idea that fixed-$\X$ theories condition on
observed regressor observations $\Xvecone$, ..., $\XvecN$, collected
in the data frame $\X$, and define a target of estimation by assuming
that the population of $\Yvec$-values at each $\Xveci$ is known:
$\Yvec_i | \Xvec_i \!\sim\! \PofYgivenXi$.  The joint distribution is
then effectively $\PofYgivenX \otimes \PhatX$, amounting to partial
plug-in of $\PhatX$ for $\PofX$ in
$\PofYandX = \PofYgivenX \otimes \PofX$.  The conditional parameter
for $\Btheta(\cdot)$ is therefore defined as
$\BthetaX = \Btheta(\PofYgivenX \otimes \PhatX)$.  We summarize
notations with emphasis on the centerline of the following box:
\begin{equation*} \label{eq:conditional-parameter}
  \arraycolsep=0pt
  \def\arraystretch{1.5}
  \begin{array}{|lcll|}
    \hline
    ~~~  \BthetaP             &~=~& \Btheta(\PofYgivenX \otimes \PofX),~~~~
                              &
    \\
    ~~~ \BthetaX              &\;=~& \Btheta(\PofYgivenX \otimes \PhatX), ~~~~
                              & ~\PhatX = (1/N)\sum \delta_{\Xvec_i} ,~~~
    \\
    ~~~ \Bthetahat            &\;=~& \Btheta(\Phat).
                              &
    \\
    \hline
  \end{array}
\end{equation*}
Note that $\X$ and $\PhatX$ contain the same information; the
conditional response distribution $\PofYgivenX$ is implied and not
shown in $\BthetaX$.  The main points are:
\begin{itemize}
\item In model-trusting theories that condition on $\X$, the target of
  estimation is~$\BthetaX$.  They assume $\BthetaX$ is the same
  for all acceptable~$\X$.
\item In model-robust theories that do not condition on $\X$, the
  target of estimation is $\Btheta(\P)$, whereas $\BthetaX$ is a
  random quantity (Corollary \ref{sec:EO} below).
\end{itemize}
The above definitions can be made more concrete by illustrating them
with the specific ways of defining regression functionals of
Section~\ref{sec:targets}:
\begin{itemize}

\item Functionals defined through minimization of objective functions:
  \begin{equation*} \label{eq:conditional-parameter-MIN}
    \arraycolsep=0pt
    \def\arraystretch{1.5}
    \begin{array}{|lcl|}
      \hline
      ~~~ \BthetaP
                        &\;=\;& \argmin_\Btheta ~ \EP[\, \cL(\Btheta;\Yvec,\Xvec) \,] ,
      \\
      ~~~ \BthetaX
                        &\;=\;& \argmin_\Btheta ~ \frac{1}{N} \sum_i \EP[\, \cL(\Btheta;\Yveci,\Xveci)\,|\,\Xveci] ,~~~
      \\
      ~~~ \Bthetahat
                        &\;=\;& \argmin_\Btheta ~ \frac{1}{N} \sum_i \cL(\Btheta;\Yveci,\Xveci) .
      \\
      \hline
    \end{array}
  \end{equation*}

\item Functionals defined through estimating equations:
  \begin{equation*} \label{eq:conditional-parameter-EE}
    \arraycolsep=0pt
    \def\arraystretch{1.5}
    \begin{array}{|llcc|}
      \hline
      ~~~ \BthetaP:~~    & \EP[\,\Bpsi(\Btheta\;\Yvec,\Xvec) \,]
                        &\;=\;& \0,~~
      \\
      ~~~ \BthetaX:~~    & \frac{1}{N} \sum_i \EP[\,\Bpsi(\Btheta;\Yveci,\Xveci)\,|\,\Xveci]
                        &\;=\;& \0,~~~
      \\
      ~~~ \Bthetahat:~   & \frac{1}{N} \sum_i \Bpsi(\Btheta;\Yveci,\Xveci)
                        &\;=\;& \0.~~
      \\
      \hline
    \end{array}
  \end{equation*}
These specialize to normal equations for linear OLS
by~\eqref{eq:OLS-moments}.
\end{itemize}
Summary: Among the three cases in each box, the most impenetrable but
also most critical case is the second one.  It defines the
``conditional parameter'' through partial plug-in of the empirical
regressor distribution.  The conditional parameter is {\em the target
  of fixed-$\X$ regression for arbitrary types of regression
  functionals}.


\subsection{Estimation Offsets}
\label{sec:EO}

The conditional parameter $\Btheta(\X)$ enables us to distinguish
between two sources of estimation uncertainty: (1)~the conditional
response distribution and (2)~the marginal regressor distribution.  To
this end we defined in Part~I for linear OLS what we call ``estimation
offsets.''  With the availability of $\Btheta(\X)$ for regression
functionals, these can be defined in full generality:
\begin{equation*}
  \arraycolsep=2pt
  \def\arraystretch{1.7}
  \begin{array}{|lcl|}
    \hline
    ~ \textit{Total~EO}
    &=&   \Bthetahat  \!-\! \BthetaP
           , ~ \\
    ~ \textit{Noise~EO}
    &=&   \Bthetahat  \!-\! \BthetaX
          , ~ \\
    ~ \textit{Approximation~EO}
    &=&   \BthetaX    \!-\! \BthetaP
          . ~ \\
    \hline
  \end{array}
\end{equation*}
The total EO is the offset of the plug-in estimate from its population
target.  The noise EO is the component of the total EO that is due to
the conditional distribution $\Yvec|\Xvec$.  The approximation EO is the
part due to the randomness of $\Xvec$ under misspecification.  These
interpretations will be elaborated in what follows.

\smallskip

\noindent{\bf Remark:} We repeat an observation made in Part~I, end of
Section~5.  The approximation EO $\BthetaX-\BthetaP$ could be
misinterpreted as a bias because it is the difference of two targets
of estimation.  This interpretation is  {\em wrong}.  In
the presence of misspecification, the approximation EO is a
non-vanishing random variable.  It will be shown to contribute not a
bias to $\Bthetahat$ but a $N^{-1/2}$ term to the sampling variability
of~$\Bthetahat$.


\subsection{Well-Specification in Terms of Approximation EOs}
\label{sec:well-spec-EO}

The approximation EO lends itself for another characterization of
well-specification:

\medskip

\noindent{\bf Proposition \ref{sec:well-spec-EO}:} {\em Assume
  $\PofX \!\mapsto \Btheta(\PofYgivenX \otimes \PofX)$ is continuous
  in the weak topology.  Then $\Bthetadot$ is well-specified for
  $\PofYgivenX\,$ iff ~$\BthetaX \!-\! \BthetaP \!=\! \0$ for all
  acceptable~$\X$.}

\medskip

\noindent {\bf Proof}: If $\Bthetadot$ is well-specified in the sense
of Section~\ref{sec:well-specification}, then
\[
\Btheta(\X) ~=~ \Btheta(\PofYgivenX \otimes \PhatX) ~=~ \Btheta(\PofYgivenX \otimes \PofX) ~=~ \Btheta(\P) .
\]
The converse follows because the empirical regressor distributions
$\PhatX$ (for $N\rightarrow\infty)$ form a weakly dense subset in the
set of all regressor distributions, and the regression functional is
assumed continuous in this argument.  ~~$\square$

A fine point about this proposition is that $\X$ is not meant as
random but as a variable taking on all acceptable regressor datasets
of arbitrarily large sample sizes.  On the other hand, here are two
consequences when $\X$ is random:

\medskip

\noindent{\bf Corollary \ref{sec:well-spec-EO}:} {\em Same assumptions
  as in Proposition~\ref{sec:well-spec-EO}.

  \begin{itemize} \itemsep 0.5em

  \item Fixed-$\X$ and random-$\X$ theories estimate the same target
    iff $\Bthetadot$ is well-specified for $\PofYgivenX$.

  \item $\Bthetadot$ is well-specified for $\PofYgivenX$ iff $ $
    $\VP[\BthetaX] = \0$ for all acceptable~$\PofX$.

  \end{itemize}

}

The first bullet confirms that the notion of well-specification for
regression functionals hits exactly the point of agreement between
theories that condition on the regressors and those that treat them as
random.  The second bullet leads the way to the fact that a
misspecified regression functional will incur sampling variability
originating from the randomness of the regressors.


\subsection{Deterministic Association Annihilates the Noise EO}
\label{sec:deterministic}

While well-specification addresses a vanishing approximation EO, one
can also consider the dual concept of a vanishing noise EO.  Here is a
sufficient condition under which the noise EO vanishes for all
regression functionals:

\medskip

\noindent {\bf Proposition \ref{sec:deterministic}:} ~ {\em If
  $\,\Yvec \!=\! f(\Xvec)\,$ is a deterministic function of $\Xvec$, then \\
  $~\Bthetahat - \BthetaX = \0$ for all regression
  functionals.}

\medskip

\noindent {\bf Proof}: The conditional response distribution is
$\PofYgivenx = \delta_{\yvec=f(\xvec)}$, hence the joint distribution
formed from $\PofYgivenx$ and $\PhatX$ is $\Phat$:
$~\PofYgivenX \otimes \PhatX = \Phat$.  It follows that
$\Btheta(\X) = \Btheta(\PofYgivenX \otimes \PhatX ) = \Btheta(\Phat) =
\Bthetahat$.~~~$\square$

\smallskip

The proposition illustrates the fact that the noise EO is due to
``noise'', that is, variability of $\Yvec$ conditional on $\Xvec$.  Thus,
although less transparent than in linear OLS, the conditional response
distribution $\Yvec|\Xvec$ is the driver of the noise~EO.




\subsection{Well-Specification and Influence Functions}
\label{sec:IC}

This section introduces influence functions for regression functionals
which will prove useful for approximations in Section~\ref{sec:EO-IC}
and for asymptotic decompositions in Section~\ref{sec:CLT}.  For
background on influence functions see, for example, Hampel et
al.~(1986)\nocite{ref:HRRS-1986} and
Rieder~(1994)\nocite{ref:Rieder-1994}.

The influence function is a form of derivative on the space of
probability distributions, which makes it an intuitive tool to
characterize well-specification of regression functionals: If
$\Btheta(\PofYgivenX \otimes \PofX)$ is constant in the argument
$\PofX$ at a fixed $\PofYgivenX$, then this means intuitively that the
``partial derivative'' wrt $\PofX$ vanishes.

The definition of the full influence function of $\Btheta(\cdot)$ is
as follows:
\begin{equation} \label{eq:IC}
\IC(\yvec,\xvec)
~=~
\frac{d}{d t} \bigg|_{t=0} \, \Btheta \left((1\!-\!t)\P + t \delta_{(\yvec,\xvec)} \right) .
\end{equation}
We omit $\Bthetadot$ as well as $\P\!=\!\PofYandX$ as arguments of
$\IC(\yvec,\xvec)$ because both will be clear from the context, except
for one occasion in Appendix~\ref{sec:partial-influence} where we
write $\IC(\yvec,\xvec;\P)$.  More relevant is the following
definition of the partial influence function of $\Btheta(\cdot)$ with
regard to the regressor distribution:
\begin{equation} \label{eq:IC-X}
\IC(\xvec)
~=~
\frac{d}{d t} \bigg|_{t=0} \, \Btheta \big(\PofYgivenX \otimes ((1\!-\!t)\PofX + t \delta_{\xvec}) \big) .
\end{equation}
For derivations of the following Lemma and Proposition, see
Appendix~\ref{sec:partial-influence}.

\medskip

\noindent{\bf Lemma \ref{sec:IC}:} {\em
$~~~\IC(\xvec) ~=~ \EP\big[\, \IC(\Yvec,\Xvec) \,|\, \Xvec\!=\!\xvec \,\big]$.
}

\medskip

\noindent{\bf Proposition \ref{sec:IC}:} {\em A regression functional
  $\Bthetadot$ with an influence function at $\PofYandX$ is
  well-specified for $\PofYgivenX$ iff
  $ ~~\IC(\xvec) = \0 ~~\forall \xvec $.  }

\medskip

\smallskip


\subsection{Approximating Estimation Offsets with Influence Functions}
\label{sec:EO-IC}

For linear OLS, Definition and Lemma~5 in Part~I exhibited an
intuitive correspondence between the total, noise and approximation EO
on the one hand and the population residual, the noise and the
nonlinearity on the other hand.  No such direct correspondence exists
for general types of regression.  The closest general statement about
EOs is in terms of approximations based on influence functions.
Assuming asymptotic linearity of $\Bthetadot$, the EOs have the
following approximations to order~$o_P(N^{-1/2})$:
\begin{equation} \label{eq:Theta-EOs}
  \arraycolsep=2pt
  \def\arraystretch{1.7}
  \begin{array}{|llcl|}
    \hline
    ~ \textit{Total~EO:}~~~~~
    &\Bthetahat  \!-\! \BthetaP
    &\approx&    \frac{1}{N} \sum_i \IC(\Yveci,\Xveci)
           , \\
    ~ \textit{Noise~EO:}
    &\Bthetahat  \!-\! \BthetaX
    &\approx&    \frac{1}{N} \sum_i \left( \IC(\Yveci,\Xveci) \!-\! \EP[\IC(\Yvec,\Xveci)|\Xveci] \right)
          , \\
    ~ \textit{Approx.~EO:}
    &\BthetaX    \!-\! \BthetaP
    &\approx&    \frac{1}{N} \sum_i \EP[\IC(\Yvec,\Xveci)|\Xveci]
          . \\
    \hline
  \end{array}
\end{equation}
These approximations \eqref{eq:Theta-EOs} lead straight to the CLTs of
the next section.


\section{Model-Robust Central Limit Theorems Decomposed}
\label{sec:CLT}

\subsection{CLT Decompositions Based on Influence Functions}
\label{sec:CLT-IF}

As in Section \ref{sec:EO-IC} assume the regression functional
$\Btheta(\P)$ is asymptotically linear with influence function
$\IC(\yvec,\xvec)$ and partial influence function
$\IC(\xvec) \!=\! \EP[\IC(\Yvec,\Xvec) \,|\Xvec\!=\!\xvec]$.  The EOs
obey the following CLTs:

\begin{equation*} \label{eq:CLT-IC}
  \arraycolsep=5pt
  \def\arraystretch{1.9}
  \begin{array}{|lcl|}
    \hline
    \sqrt{N} \, (\Bthetahat - \BthetaP) \!\!
    &\stackrel{\Dist}{\longrightarrow}&
    \Norm\left(\zero, \, \VP[\IC(\Yvec,\Xvec)] \right) ,
    \\
    \sqrt{N} \, (\Bthetahat - \Btheta(\X)) \!\!
    &\stackrel{\Dist}{\longrightarrow}&
    \Norm\left(\zero, \, \EP[\VP[\IC(\Yvec,\Xvec) \,|\,\Xvec]] \right) ,
    \\
    \sqrt{N} \, (\Btheta(\X) - \BthetaP) \!\!
    &\stackrel{\Dist}{\longrightarrow}&
    \Norm\left(\zero, \, \VP[\EP[\IC(\Yvec,\Xvec) \,|\,\Xvec]] \right) .
    \\
    \hline
    \end{array}
\end{equation*}

\smallskip

\noindent
These are immediate consequences of the assumed asymptotic linearities.
The asymptotic variances of the EOs follow the canonical decomposition
\begin{equation*}
  \VP[\IC(\Yvec,\Xvec)]
  ~=~
  \EP[\VP[\IC(\Yvec,\Xvec)\,|\,\Xvec]\,]
  +
  \VP[\EP[\IC(\Yvec,\Xvec)\,|\,\Xvec]\,] ,
\end{equation*}
the three terms being the asymptotic variance-covariance matrices of
the total, the noise and the approximation EO, respectively.  Implicit
in this Pythagorean formula is that
$\IC(\Yvec,\Xvec) \!-\! \EP[\IC(\Yvec,|\Xvec)$ and
$\EP[\IC(\Yvec,|\Xvec)$ are orthogonal to each other, which implies by
\eqref{eq:Theta-EOs} that the noise EO and the approximation EO are
asymptotically orthogonal.  Asymptotic orthogonalities based on
conditioning are well-known in semi-parametric theory.  For linear OLS
this orthogonality holds exactly for finite $N$ due to (6) and (13) in
Part~I:
$~\VP[\Bbetahat \!-\! \Bbeta(\X), \Bbeta(\X) \!-\! \Bbeta(\P)] = \0$.

The following corollary is a restatement of Proposition~\ref{sec:IC},
but enlightened by the fact that it relies on the asymptotic variance of
the approximation~EO.

\medskip

\noindent{\bf Corollary \ref{sec:CLT-IF}:} {\em The regression
  functional $\Bthetadot$ is well-specified for $\PofYgivenX$ iff the
  asymptotic variance of the approximation EO vanishes for all
  acceptable~$\PofX$.}

\medskip

\noindent {\bf Proof}: Using careful notation the condition says
$\V_\PofX[\E_\PofYgivenX[\IC(\Yvec,\Xvec)|\Xvec]] = \0$ for all
acceptable~$\PofX$.  This in turn means
$\E_\PofYgivenX[\IC(\Yvec,\Xvec)|\Xvec=\xvec] = \0$ for all~$\xvec$, which
is the condition of Proposition~\ref{sec:IC}. ~~$\square$

\subsection{CLT Decompositions for EE Functionals}

For EE functionals the influence function is
$\IC(\yvec,\xvec) = \BLambda(\Btheta)^{-1} \Bpsi(\Btheta;\yvec,\xvec)$ where
$\Btheta = \Btheta(\P)$ and
$\BLambda(\Btheta) = \nabla_\Btheta \EP[\Bpsi(\Btheta;\Yvec,\Xvec)]$
is the Jacobian of size $q\times q$,
$q\!=\!\dim(\Bpsi)\!=\!\dim(\Btheta)$.  Then the CLTs specialize to
the following:
\begin{equation*} \label{eq:Huber-CLT}
  \arraycolsep=5pt
  \def\arraystretch{1.9}
  \begin{array}{|rcl|}
    \hline
    \sqrt{N} \, (\Bthetahat        - \Btheta            ) \!\!
    &\stackrel{\Dist}{\longrightarrow}& \!\!
    \Norm\left(\zero, \, \BLambda(\Btheta)^{-1} \, \VP[\Bpsi(\Btheta;\Yvec,\Xvec)] \; \BLambda(\Btheta)\Tr^{-1} \right)
    ~\\~~
    \sqrt{N} \, (\Bthetahat        - \BthetaX) \!\!
    &\stackrel{\Dist}{\longrightarrow}& \!\!
    \Norm\left(\zero, \, \BLambda(\Btheta)^{-1} \, \EP[\VP[\Bpsi(\Btheta;\Yvec,\Xvec)\,|\,\Xvec]\,] \; \BLambda(\Btheta)\Tr^{-1} \right)
    ~\\~~
    \sqrt{N} \, (\BthetaX - \Btheta     ) \!\!
    &\stackrel{\Dist}{\longrightarrow}& \!\!
    \Norm\left(\zero, \, \BLambda(\Btheta)^{-1} \, \VP[\EP[\Bpsi(\Btheta;\Yvec,\Xvec)\,|\,\Xvec]\,] \; \BLambda(\Btheta)\Tr^{-1} \right)
    ~\\
    \hline
  \end{array}
\end{equation*}

\smallskip

\noindent The first line is Huber's (1967,
Section~3)\nocite{ref:Huber-1967} result.  The asymptotic variances
have the characteristic sandwich form.  It is natural that they are
related according to
\begin{equation*}
  \VP[\Bpsi(\Btheta;\Yvec,\Xvec)]
  ~=~
  \EP[\VP[\Bpsi(\Btheta;\Yvec,\Xvec)\,|\,\Xvec]\,]
  +
  \VP[\EP[\Bpsi(\Btheta;\Yvec,\Xvec)\,|\,\Xvec]\,] ,
\end{equation*}
where on the right side the first term relates to the noise EO and the
second term to the approximation~EO.

Linear OLS is a special case with
$~\Bpsi(\Bbeta;y,\xvec) \!=\!\xvec \xvec\Tr \Bbeta - \xvec y$,
$\BLambda \!=\! \EP[\Xvec \Xvec\Tr]$,
$\IC(y,\xvec) \!=\! \EP[\Xvec \Xvec\Tr]^{-1} (\xvec \xvec\Tr \Bbeta - \xvec
y)$, and hence the CLTs of Part~I, Proposition~7.


\subsection{Implications of the CLT Decompositions}
\label{sec:model-bias}

We address once again potential confusions relating to different
notions of bias.  Misspecification, in traditional parametric
modeling, is sometimes called ``model bias'' which, due to unfortunate
terminology, may suggest a connection to estimation bias,
$\EP[\Bthetahat_N] - \Btheta(\P)$.  Importantly, there is no
connection between the two notions of bias.  Estimation bias typically
vanishes at a rate faster than $N^{-1/2}$ and does not contribute to
standard errors derived from asymptotic variances.  Model bias, on the
other hand, which is misspecification, generates in conjunction with
the randomness of the regressors a contribution to the standard error,
and this contribution is asymptotically of order $N^{-1/2}$, the same
order as the better known contribution due to the conditional noise in
the response.  This is what the CLT decomposition shows.  It also
shows that the two sources of sampling variability are asymptotically
orthogonal.  ---~In summary:

\begin{itemize}
\item[] {\em Model bias/misspecification does not create estimation
    bias; it creates sampling variability to the same order as the
    conditional noise in the response. }
\end{itemize}

\smallskip


\section{Plug-In/Sandwich Estimators versus $M$-of-$N$ Bootstrap Estimators of Standard Error}
\label{sec:plugin-bootstrap}

\subsection{Plug-In Estimators are Limits of $M$-of-$N$ Bootstrap Estimators }
\label{sec:plugin-as-limit}

In Part~I, Section~8, it was indicated that for linear OLS there
exists a connection between two ways of estimating asymptotic
variance: the sandwich estimator for sample size $N$ is the limit of
the $M$-of-$N$ bootstrap as $M \rightarrow \infty$, where bootstrap is
the kind that resamples $\xy$ cases rather than residuals.  This
connection holds at a general level: all plug-in estimators of
standard error are limits of bootstrap in this sense.

The crucial observation of Part~I goes through as follows: The
$M$-of-$N$ bootstrap is iid sampling of $M$ observations from some
distribution, hence there must hold a CLT as the resample size grows,
$M \rightarrow \infty$.  The distribution being (re)sampled is the
empirical distribution $\PhatN = (1/N) \sum \delta_{(\yi,\xveci)}$,
where $N$ is fixed but $M \rightarrow \infty$.\footnote{This causes
  ever more ties at $M$ grows.}  Therefore, the following holds for
bootstrap resampling of any well-behaved statistical functional, be it
in a regression context or not:

\medskip

\noindent{\bf Proposition~\ref{sec:plugin-as-limit}:} {\em Assume
  the regression functional $\Bthetadot$ is asymptotically normal
  for a sufficiently rich class of joint distributions~$\P=\PofYandX$ with
  acceptable regressor distributions~$\PofX$ as follows:
  \[
  N^{1/2} (\BthetahatN - \Btheta(\P))
    ~~\overset{\Dist}{\longrightarrow}~~
    \Norm\left(\zero,~ \AV[\P;\Bthetadot] \right)
    ~~~~(N\rightarrow\infty).
  \]
  Let a fixed dataset of size $N$ with acceptable regressors be
  represented by the empirical measure $\PhatN$.
  Then a CLT holds for the $M$-of-$N$ bootstrap as $M \rightarrow \infty$,
  with an asymptotic variance obtained by plug-in.  Letting
  $\BthetabootM = \Btheta(\PbootM)$ where $\PbootM$ is the empirical distribution
  of a resample of size $M$ from $\PhatN$, we have:}
  \begin{equation*} \label{eq:empirical-CLT}
    M^{1/2} \, (\BthetabootM - \BthetahatN)
    ~~\overset{\Dist}{\longrightarrow}~~
    \Norm\left(\zero,~ \AV[\PhatN;\Bthetadot] \right)
    ~~~~(M\rightarrow\infty,~~N~\textit{fixed}).
\end{equation*}

\noindent
The proposition contains its own proof.  The following is the
specialization to EE functionals where the asymptotic variance has the
sandwich form:

\medskip

\noindent{\bf Corollary~\ref{sec:plugin-as-limit}:}{\em ~The plug-in
  sandwich estimator for an EE functional is the asymptotic variance
  estimated by the $M$-of-$N$ bootstrap in the limit
  $M \!\!\rightarrow\!\infty$ for a fixed sample of size~$N$.  }


\subsection{Arguments in Favor of  $M$-of-$N$ Bootstrap Over Plug-In Estimators}

A natural next question is whether the plug-in/sandwich estimator is
to be preferred over $M$-of-$N$ bootstrap estimators, or whether there
is a reason to prefer some form of $M$-of-$N$ bootstrap.  In the
latter case the follow-up question would be how to choose the resample
size~$M$.  While we do not have any recommendations for choosing a
specific~$M$, there exist various arguments in favor of some
$M$-of-$N$ bootstrap over plug-in/sandwich estimation of standard error.

A first argument is that bootstrap is more flexible in that it lends
itself to various forms of confidence interval construction that grant
higher order accuracy of coverage.  See, for example, Efron and
Tibshirani (1994)\nocite{ref:ET-1994} and Hall
(1992)\nocite{ref:Hall-1992}.

A second argument is related to the first but in a different
direction: Bootstrap can be used to diagnose whether the sampling
distribution of a particular functional $\Bthetadot$ is anywhere near
asymptotic normality for a given sample size~$N$.  This can be done by
applying normality tests to simulated bootstrap values $\Bthetabootb$
$(b=1,...,B)$, or by displaying these values in a normal quantile
plot.

A third argument is that there exists theory that shows bootstrap to
work for very small $M$ compared to $N$ in some situations where even
conventional $N$-of-$N$ bootstrap does not work.  (See Bickel, G\"otze
and van Zwet (1997)\nocite{ref:BGvZ-1997} following Politis and Romano
(1994)\nocite{ref:PR-1994} on subsampling.)  It seems therefore
unlikely that the limit $M\!\rightarrow\!\infty$ for fixed $N$ will
yield any form of superiority to bootstrap with finite~$M$.

A fourth argument derives from a result by Buja and Stuetzle
(2016)\nocite{ref:BS-2016}, which states that so-called ``$M$-bagged
functionals'' have low complexity in a certain sense, the lower the
smaller the resample size $M$ is.  The limit
$M \!\rightarrow\! \infty$ is therefore the most complex choice.  The
connection to the issue of ``bootstrap versus plug-in/sandwich
estimators'' is that $M$-of-$N$ bootstrap standard errors are simple
functions of $M$-bagged functionals, hence the complexity comparison
carries over to standard errors.

It appears that multiple arguments converge on the conclusion that the
$M$-of-$N$ bootstrap is to be preferred over plug-in/sandwich standard
errors.  Also recall that both are to be preferred over the residual
bootstrap.


\section{Summary and Conclusion}
\label{sec:summary}

This article completes important aspects of the program set out in
Part~I.  It pursues the idea of model robustness to its conclusion for
arbitrary types of regression based on iid observations.  The notion
of model robustness coalesces into a model-free theory where all
quantities of interest are statistical functionals, called
``regression functionals'', and models take on the role of heuristics
to suggest objective functions whose minima define regression
functionals defined on largely arbitrary joint $(\Yvec,\Xvec)$
distributions.  In this final section we recount the path that makes
the definition of well-specification for regression functionals
compelling.

To start, an important task of the present article has been to extend
the two main findings of Part~I from linear OLS to arbitrary types of
regression.  The findings are that nonlinearity and randomness of the
regressors interact (``conspire'')
\begin{itemize}
\item[(1)] to cause the target of estimation to depend on the regressor
  distribution;
\item[(2)] to cause $N^{-1/2}$ sampling variability to arise that is
  wholly different from the sampling variability caused by the
  conditional noise in the response.
\end{itemize}
It was intuitively clear that these effects would somehow carry over
from linear OLS to all types of regression, but it wasn't clear what
would take the place of ``nonlinearity,'' a notion of first order
misspecification peculiar to fitting linear equations and estimating
linear slopes.  In attempting to generalize Part~I, a vexing issue is
that one is looking for a framework free of specifics of fitted
equations and additive stochastic components of the response.
Attempts at directly generalizing the notions of ``nonlinearity'' and
``noise'' of Part~I lead to dead ends of unsatisfactory extensions
that are barely more general than linear OLS.  This raises the
question to a level of generality in which there is very little air to
breathe: the objects that remain are a regression functional
$\Bthetadot$ and a joint distribution~$\PofYandX$.  Given these two
objects, what do mis- and well-specification mean?  An answer, maybe
{\em the} answer, is arrived at by casting regression in the most
fundamental way possible: {\em Regression is the attempt to describe
  the conditional response distribution~$\PofYgivenX$.}  This
interpretation sweeps away idiosyncratic structure of special cases.
It also suggests taking the joint distribution $\PofYandX$ apart and
analyzing the issue of mis- and well-specification in terms of
$\PofYgivenX$ and~$\PofX$, as well as $\Bthetadot$, the quantities of
interest.  The solution, finally, to
\begin{itemize}
\item establishing a compelling notion of mis- and well-specification
  at this level of generality, and
\item extending (1) and (2) above to arbitrary types of regression,
\end{itemize}
is to look no further and use the ``conspiracy effect'' (1) as the
definition: Misspecification means dependence of the regression
functional on the regressor distribution.  Conversely,
well-specification means the regression functional does not depend on
the regressor distribution; it is a property of the conditional
response distribution alone.

The ``conspiracy effect'' (2) above is now a corollary of the
definition: If the functional is not constant across regressor
distributions, it will incur random variability on empirical regressor
distributions, and this at the familiar rate~$N^{-1/2}$.

The link between the proposed definition and conventional ideas of
mis-/well-specification is as follows: Because most regressions
consist of fitting some functional form of the regressors to the
response, misspecification of the functional form is equivalent to
misspecification of its parameters viewed as regression functionals:
depending on where the regressors fall, the misspecified functional
form needs adjustment of its parameters to achieve the best
approximation over the distribution of the regressors.

Well-specification being an ideal, in reality we always face degrees
of misspecification.  Acknowledging the universality of
misspecification, however, does {\em not} justify carelessness in
practice.  It is mandatory to perform diagnostics and, in fact, we
proposed a type of diagnostic in Sections~\ref{sec:reweighting} and
\ref{sec:reweighting-methodology} tailored to the present notion of
mis-/well-specification.  The diagnostic consists of checking the
dependence of regression functionals on the regressor distribution by
systematically perturbing the latter, not by shifting or otherwise
moving it, but by reweighting it.  Reweighting has the considerable
advantage over other forms of perturbation that it applies to all
variable types, not just quantitative ones.

While the reality of misspecification imposes a duty to perform
diagnostics, there is also an argument to be made to feel less guilty
about choosing simpler models over more complex ones.  One reason is
that the reweighting diagnostic permits localization of models and
thereby enables a systematic exploration of local best approximations,
always in terms of model parameters interpreted as regression
functionals.  As shown in
Sections~\ref{sec:focal-coefficient}-\ref{sec:focal-reweighting}, this
possibility vastly extends the expressive power of models beyond that
of a single model fit.

Finally, there is an argument to be made in favor of using statistical
inference that is model-robust, and to this end one can use $\xy$
bootstrap estimators or plug-in/sandwich estimators of standard
errors.  Between the two, one can give arguments in favor of bootstrap
over plug-in/sandwich estimators.  Most importantly, though, both
approaches to inference are in accord with the insight that
misspecification forces us to treat regressors as random.


\smallskip

\noindent{\bf Acknowledgments}: We are grateful to Gemma Moran and
Bikram Karmakar for their help in the generalizations of
Section~\ref{sec:conditional-parameter}, and to Hannes Leeb for
pointing out the source of the Tukey quote shown before the
introduction.



\newpage


\appendix


\section*{Appendix}

\subsection{Assumptions}
\label{sec:assumptions}

When defining a regression functional $\Btheta(\P)$, one needs to
specify a set $\cP$ of joint distributions $\P \!=\! \PofYandX$ for
which the functional is defined.  This set can be specific to the
functional in several ways.  Here is a list of conditions on $\cP$
that will be assumed as needed:
\begin{itemize}
\item[a)] Expectations $\EP[f(\Xvec,\Yvec)]$ exist as needed for all
  distributions in~$\P \in \cP$.
\item[b)] If the regression functional derives from parameters of
  fitted equations, it will be assumed that the regressor distribution
  $\PofX$ grants identifiabiliy of the fitted parameters, as when
  strict collinearity of the regressor distribution needs to be
  excluded in order to uniquely fit linear equations.  If this is the
  case we will say $\PofX$ is an ``{\em acceptable}'' regressor
  distribution.
\item[c)] With b) in mind, we will assume that if a regressor
  distribution $\PofX$ is acceptable, a mixture
  $\alpha \PofX + (1\!-\!\alpha) \PofXprime$ with any other
  distribution $\PofXprime$, subject to a) above, will also be
  acceptable, the reason being that mixing can only enlarge but not
  diminish the support of the distribution, hence identifiability will
  be inherited from~$\PofX$ irrespective of~$\PofXprime$.
\item[d)] Item c) ensures that the set of acceptable regressor
  distributions is so rich that $E_\PofX[f(\Xvec] = 0$ for all
  acceptable $\PofX$ entails $f \equiv 0$. Reason: Mix with atoms at
  arbitrary locations.
\item[e)] For working models
  $\{ \Qmodel\!\!: \Btheta \!\in\! \BTheta \}$ (not treated as
  correct) it will be assumed $\Qmodel \otimes \PofX \in \cP$ for all
  $\Btheta \!\in\! \BTheta$ and all acceptable~$\PofX$.
\item[f)] Where conditional model densities $q(\yvec|\xvec;\Btheta)$ of
  the response appear, they will be densities with regard to some
  dominating measure~$\nu(d\yvec)$.
\item[g)] For plug-in estimation it will be required that for $N$ iid
  draws $(\Yveci,\Xveci) \sim \P \!\in\! \cP$ the empirical distribution
  \[
  \Phat = \PhatN = (1/N) \sum \delta_{(\Yveci,\Xveci)}
  \]
  is in $\cP$ with limiting probability~1 as
  $N\!\!\rightarrow\!\!\infty$.  For example, one needs
  $N\!\ge\!p\!+\!1$ non-collinear observations in order to fit a
  linear equation with intercept.
\item[h)] To form influence functions for regression functionals, it will
  be assumed that for $\P \!\in\! \cP$ and
  $(\yvec,\xvec) \!\in\! \cX \!\times\! \cY$ we have
  $(1\!-\!t) \P + t \delta_{y,\xvec} \in \cP$ for $0\!<\!t\!<\!1$.
\item[i)] $\cP$ will be assumed to be convex, hence closed under finite
  mixtures.
\end{itemize}

\subsection{Proper Scoring Rules, Bregman Divergences and Entropies}
\label{ref:proper-scoring}

\subsubsection{General Theory:}

We describe objective functions called ``proper scoring rules'', which
generalize negative log-likelihoods, and associated Bregman
divergences, which generalize Kullback-Leibler (K-L) divergences.
Proper scoring rules can be used to extend the universe of regression
functionals based on working models.  For insightful background on
proper scoring rules, see Gneiting and Raftery
(2007)\nocite{ref:GR-2007} (but note two reversals of conventions:
theirs are gain functions where ours are loss functions, and their
order of arguments is switched from ours).

We begin by devising discrepancy measures between pairs of
distributions based on axiomatic requirements that can be gleaned from
two properties of K-L divergences: Fisher consistency at the working
model and availability of plug-in for ``empirical risk minimization''
(machine learning terminology).  The resulting functionals will be
called ``proper scoring functionals.''

Denoting a discrepancy measure between two distributions by
$D(\P,\Q)$, the intended roles of the two arguments are that $\P$ is
the actual data distribution and $\Q$ is a member of a working model.
Fisher consistency of a minimum discrepancy functional follows from
the following requirements:
\begin{itemize}
\item[] (A)~~$D(\P,\Q) \!\ge\! 0$ with equality iff $\P\!=\!\Q$.
\end{itemize}
The ``only if'' part in the second clause is essential.  Other
properties such as symmetry and triangle inequalities are not needed
and would be too restrictive.

As for the availability of plug-in estimation, it would follow from a
structural property such as dependence of $D(\P,\Q)$ on $\P$ only
through its expectation $\EP[...]$, whose plug-in estimate is the
empirical mean.  Other types of plug-in exist, for example for
quantiles, in particular medians.  Yet other cases, such as Hellinger
distances, require density estimation for plug-in, which adds a layer
of complexity.  In what follows we will impose the strong condition
that $D(\P,\Q)$ depends on $\P$ essentially only through $\EP[\cdot]$,
but this requirement only concerns the part of $D(\P,\Q)$ that is
relevant for minimization over working model distributions~$\Q$.  We
can use the K-L divergence as a guide: In
\[
\KL(\P,\Q) ~=~ \EP[-\!\log \frac{q(\Yvec)}{p(\Yvec)}] ~=~ \EP[-\!\log q(\Yvec)] - \EP[-\!\log p(\Yvec)],
\]
the second term requires for plug-in a density estimate of $p(\yvec)$,
but this term does not depend on $\Q$, hence is irrelevant for
minimization over~$\Q$.  By analogy we impose the following structural
form on the discrepancy measure:
\begin{itemize}
\item[] (B)~~$D(\P,\Q) ~=~ \EP[S(\Yvec,\Q)] - H(\P)$ .
\end{itemize}
This condition, combined with condition (A), constrains $D(\P,\Q)$ to
be a so-called ``{\bf Bregman divergence}''.  The following structure
falls into place:
\begin{itemize} \itemsep 0.5em
\item Define $S(\P,\Q) = \EP[S(\Yvec,\Q)]$.  Then $S(\P,\P) = H(\P)$
  due to~(A).
\item The term $S(\Yvec,\Q)$ is a so-called ``{\bf strict proper scoring
    rule}'', characterized by
  $S(\P,\Q) \ge S(\P,\P)$, with equality iff $\P = \Q$.  This is a
  direct translation of (A) applied to
  $D(\P,\Q) = S(\P,\Q) - S(\P,\P)$.
\item The term $H(\P)$ is an ``{\bf entropy}'' as it is a strictly
  concave functional of~$\P$.  Its upper tangent at tangent point $\Q$
  is $\P \mapsto S(\P,\Q)$ due to~(A).  Also, (A) excludes tangent
  points other than $\Q$, hence renders $H(\P)$ strictly concave.
\end{itemize}
Strict proper scoring rules $S(\yvec,\Q)$ generalize negative
log-likelihoods.

\subsubsection{Examples of Proper Scoring Rules --- Density Power Divergences:}

A one-parameter family of strict proper scoring rules is as follows:
\[
S_\alpha(\yvec,\Q) =
\left\{
\begin{array}{lll}
- q^{\alpha}(\yvec)/\alpha + \int q^{1+\alpha} \, d\mu \,/(1+\alpha)  & \textrm{for} & \alpha \neq 0, -1, \\
- \log(q(\yvec))  & \textrm{for} & \alpha = 0, \\
~~ 1/q(\yvec) + \int \log(q) \, d\mu & \textrm{for} & \alpha = -1.
\end{array}
\right.
\]
These include proper scoring rules derived from the ``density power
divergences'' of Basu et al.~(1998)\nocite{ref:Basu-et-al-1998} for
$\alpha \!>\! 0$, the negative log-likelihood for $\alpha \!=\! 0$,
and a proper scoring rule derived from the Itakura-Saito divergence
for $\alpha\!=\!-1$.  The two logarithmic cases ($\alpha\!=\!0,1$)
form smooth fill-in in the manner of the logarithm in the Box-Cox
family of power transforms, which makes the family well-defined for
all $\alpha \!\in\! \Reals$.  The case $\alpha \!=\! 1$ corresponds to
the $L_2$ distance $D_2(\P,\Q) = \int (p - q)^2 d\mu$; its proper
scoring rule is $S(\yvec,\Q) = -q(\yvec) + \int q^2 \, d\mu /2$ and its
entropy is the Gini index $H(\P) = -\!\int p^2 \,d\mu /2$.  The power
$\alpha$ is a robustness parameter, in the meaning of insensitivity to
tails: robustness is gained for $\alpha \uparrow$ and sensitivity to
tail probabilities for $\alpha \downarrow$.  Basu et
al.~(1998)\nocite{ref:Basu-et-al-1998} show that for $\alpha \!>\! 0$
the influence function is redescending for the minimum divergence
estimator of the normal working model.  For $\alpha \!\le\! -1$ the
divergence is so sensitive to small probabilities (hence the opposite
of robust) that model densities $q(\yvec)$ need to have tails lighter even
than normal distributions.


\subsubsection{Proper Scoring Rules for Regression:}
\label{sec:PS-for-regression}

When applying a proper scoring rule $S(\yvec,\Q)$ to regression, scoring
is on the conditional response distributions
$\Qmodelx$ $=$ $\Q(d\yvec|\xvec;\Btheta)$ in light of a response value~$\yvec$
at~$\xvec$.  The resulting objective function is therefore:
\[
\cL(\Btheta; \yvec,\xvec) ~=~ S(\yvec, \Qmodelx) \,,
\]
which is used to construct a regression functional with argument
$\P = \PofYandX$ by
\[
\Btheta(\P) = \argmin_{\Btheta \in \BTheta} \, \EP[ \cL(\Btheta; \Yvec,\Xvec) ] \,.
\]
Fisher consistency follows from the fact that if $\PofYgivenx \!=\! \Qmodelxzero$,
then $\Btheta_0$ minimizes the objective function conditionally at each~$\xvec$
due to proper scoring:
\begin{equation} \label{eq:pointwise-min}
\E_\PofYgivenx [ \cL(\Btheta_0; \Yvec, \xvec) ] \le  \E_\PofYgivenx [ \cL(\Btheta; \Yvec, \xvec) ]
~~~ \forall \Btheta,~\forall \xvec.
\end{equation}
The same holds after averaging over arbitrary regressor
distributions~$\PofX(d\xvec)$:
\[
\EP [ \cL(\Btheta_0; \Yvec, \Xvec) ] \le  \EP [ \cL(\Btheta; \Yvec, \Xvec) ]  ~~~ \forall \Btheta,
\]
and hence $\Btheta(\P) = \Btheta_0$.

\bigskip

\subsubsection{Pointwise Bregman Divergences from Convex Functions:}

We illustrate one simple way of constructing what one may call
``pointwise'' Bregman divergences to convey the role of convex
geometry.  (We use here convex rather than concave functions, but this
is immaterial for the construction.)  If $\phi(q)$ is a strictly
convex smooth function, define the associated discrepancy between two
values $p$ and $q$ (in this order) to be
$d(p,q) = \phi(p) - (\phi(q) + \phi'(q) (p-q))$.  The term in parens
is the subtangent of $\phi(\cdot)$ at $q$ as a function of $p$, hence
$d(p,q) \!\ge\! 0$ holds due to convexity, and $d(p,q) \!=\! 0$ iff
$p\!=\!q$ due to strict convexity.  Note $d(p,q)$ is {\em not}
generally symmetric in its arguments.  The associated Bregman
divergence between distributions $\P$ and $\Q$ is obtained by applying
$d(p,q)$ to the respective densities $p(\yvec)$ and $q(\yvec)$,
integrated wrt the dominating measure $\nu(d\yvec)$:
\begin{eqnarray*}
D(\P,\Q) &=& \int \phi(p(\yvec)) \nu(d\yvec) - \int \phi(q(\yvec)) \nu(d\yvec) - \int \phi'(q(\yvec)) (p(\yvec)-q(\yvec)) \nu(d\yvec) \\
         &=&  - H(\P) + H(\Q) - \EP[\phi'(q(\Yvec))] + \E_\Q[\phi'(q(\Yvec))] ,
\end{eqnarray*}
where $H(\Q) = - \int \phi(q(\yvec)) \nu(d\yvec)$ is the associated entropy and
\[
S(\yvec,\Q) = - \phi'(q(\yvec)) + \E_\Q[\phi'(q(\Yvec))] + H(\Q).
\]
Special cases: K-L divergence for $\phi(q) \!=\! q \log(q)$; $L_2$
distance for $\phi(q) \!=\! q^2$.

\subsubsection{Density Power Divergences in Greater Detail:}

Applying the preceding subsection to power transformations, suitably
transformed to convexity following the Box-Cox transformation scheme,
one obtains the family of density power divergences.  The following is
a one-parameter family of convex functions defined {\em for all}
$\alpha\!\in\!\Reals$:
\[
\phi_\alpha(q) =
\left\{
\begin{array}{lll}
q^{1+\alpha}/(\alpha(1\!+\!\alpha)) - q/\alpha ~+~ 1/(1\!+\!\alpha) & \textrm{for} & \alpha \neq 0, -1, \\
\;q \log(q) -q + 1 & \textrm{for} & \alpha = 0, \\
- \log(q)  + q - 1 & \textrm{for} & \alpha = -1,
\end{array}
\right.
\]
The linear terms in $q$ and the constants are irrelevant but useful to
normalize $\phi_\alpha(1) \!=\! 0$ and $\phi_\alpha'(1) \!=\! 0$ for
all $\alpha \!\in\! \Reals$ and to achieve the logarithmic limits for
$\alpha\!=\!0$ and $\alpha \!=\! -1$.
The derivatives are:
\[
\phi_\alpha'(q) =
\left\{
\begin{array}{lll}
q^{\alpha}/\alpha - 1/\alpha & \textrm{for} & \alpha \neq 0, -1, \\
\;\log(q) & \textrm{for} & \alpha = 0, \\
- 1/q  + 1 & \textrm{for} & \alpha = -1,
\end{array}
\right.
\]
The associated Bregman discrepancies are:
\[
d_\alpha(p,q) =
\left\{
\begin{array}{lll}
p^{1+\alpha} / (\alpha(1\!+\!\alpha)) ~+~ q^{1+\alpha} / (1\!+\!\alpha) - p q^{\alpha} / \alpha & \textrm{for} & \alpha \neq 0, -1, \\
\;p \log(p/q) + q - p & \textrm{for} & \alpha = 0, \\
  - \log(p/q)  + p/q  & \textrm{for} & \alpha = -1,
\end{array}
\right.
\]
Integrated to form Bregman divergences for pairs of densities $p\!=\!p(\yvec)$
and $q\!=\!q(\yvec)$ of $\P$ and $\Q$, respectively, one obtains:
\[
D_\alpha(\P,\Q) =
\left\{
\begin{array}{lll}
\int \left( p^{1+\alpha}/(\alpha(1\!+\!\alpha)) ~+~ q^{1+\alpha}/(1\!+\!\alpha) - p q^{\alpha} / \alpha \right) d\mu) & \textrm{for} & \alpha \neq 0, -1, \\
\int p \log(p/q) \, d\mu & \textrm{for} & \alpha = 0, \\
- \int \left( \log(p/q)  + p/q \right) d\mu  & \textrm{for} & \alpha = -1,
\end{array}
\right.
\]
The proper scoring rules associated with density power divergences
(neglecting constants) are as follows:
\[
S_\alpha(\yvec,\Q) =
\left\{
\begin{array}{lll}
- q^{\alpha}(\yvec)/\alpha + \int q^{1+\alpha} \, d\mu \,/(1+\alpha)  & \textrm{for} & \alpha \neq 0, -1, \\
- \log(q(\yvec))  & \textrm{for} & \alpha = 0, \\
~~ 1/q(\yvec) + \int \log(q) \, d\mu & \textrm{for} & \alpha = -1.
\end{array}
\right.
\]
The associated entropies are as follows:
\[
H_\alpha(\Q) =
\left\{
\begin{array}{lll}
- \int q^{1+\alpha} \, d\mu /(\alpha(1\!+\!\alpha)) & \textrm{for} & \alpha \neq 0, -1, \\
- \int q \log(q) \, d\mu & \textrm{for} & \alpha = 0, \\
~~ \int \log(q) \, d\mu  & \textrm{for} & \alpha = -1.
\end{array}
\right.
\]


\subsection{Partial Influence Functions with regard to Regressor Distributions}
\label{sec:partial-influence}

Remark on notation: We have a need to explicitly note the distribution
at which the influence function is created.  Recall the definition
from Section~\ref{sec:IC}:
\[
\IC(\yvec,\xvec;\P) ~=~
\frac{d}{d t} \big|_{t=0} \, \Btheta((1\!-\!t)\P + t \delta_{(\yvec,\xvec)}) .
\]
This definition can be mapped to the interpretation of regression
functionals as having two separate arguments, $\PofYgivenX$ and
$\PofX$ by splitting the pointmass $\delta_{(\yvec,\xvec)}$ off to the two
arguments: The conditional response distribution is
$(1\!-\!t) \PofYgivenx + t \delta_y$ at this particular $\xvec$,
leaving those at all other $\xvec'$ unchanged; the regressor
distribution is changed to $(1\!-\!t) \PofX + t \delta_\xvec$.

We show that the partial influence functions wrt $\PofX$ is a shown in
Proposition~\ref{sec:IC}.  We start with the integrated form of the
derivative:
\[
\IC(\P';\P) ~=~
\frac{d}{d t} \big|_{t=0} \, \Btheta((1\!-\!t)\P + t \P')
~=~ \int \IC(\yvec,\xvec;\P) \, \P'(d\yvec,d\xvec) .
\]
which uses the fact that $\int \IC(\Yvec,\Xvec;\P) d\P = \0$.  To form the
partial influence function wrt $\PofX$ holding $\PofYgivenX$ fixed, we
rewrite the expansion with $\PofYgivenX$ being the same for $\P'$
and~$\P$:
\begin{equation} \label{eq:partial-IC-expantion}
\frac{d}{d t} \big|_{t=0} \, \Btheta(\PofYgivenX \otimes ((1\!-\!t)\PofX'+t\PofX))
~=~ \int \int \IC(\yvec,\xvec) \, \P(d\yvec|d\xvec) \P'(d\xvec) ,
\end{equation}
which shows that the partial influence function wrt $\PofX$ is
\[
\IC(\xvec;\PofX) ~=~ \EP[ \IC(\Yvec,\Xvec;\P) | \Xvec=\xvec ] .
\]
(We assumed that if $\PofX$ is an acceptable regressor distribution, so
is a mixture $(1\!-\!t)\PofX + t \delta_\xvec$ for small~$t>0$ and
any~$\xvec$.)

To show Proposition~\ref{sec:IC}, if we have well-specification, then
$\Btheta((1\!-\!t)\P + t \delta_\xvec) = \Btheta(\P)$, hence
$\IC(\xvec;\P) = 0$.  For the converse, we use the following integral
representation, which is integrating up derivatives along a convex
segment:
\[
\Btheta(\PofYgivenX \otimes \PofX') ~=~ \Btheta(\PofYgivenX \otimes \PofX) + \int \IC(\PofX'; (1\!-\!t)\PofX+t\PofX') \, dt .
\]
As a consequence, if $\IC(\xvec; (1\!-\!t)\PofX+t\PofX') = 0$ for all
$\xvec$ at all regressor distributions, then
$\IC(\PofX'; (1\!-\!t)\PofX+t\PofX') = 0$ for all $\PofX'$ and
$\PofX$, hence
$\Btheta(\PofYgivenX \otimes \PofX') = \Btheta(\PofYgivenX \otimes \PofX)$ for all
$\PofX'$ and~$\PofX$.  ~~$\square$

\bigskip





\subsection{Proof of Lemma \ref{sec:well-spec-ML-PS-EE}.3}
\label{sec:proof}

\noindent The ``if'' part is trivial as it involves taking
expectations wrt arbitrary $\PofX$.  The ``only if'' part follows by
observing that for any acceptable $\PofX$ with
$\Btheta(\PofYgivenX \otimes \PofX) = \Btheta_0$ there must exist
$\xvec$ for which
$\EP[ \Bpsi(\Btheta_0;\Yvec,\Xvec) | \Xvec\!=\!\xvec ] \neq \0$.
Mixtures $\PtX = (1\!-\!t)\PofX + t \delta_\xvec$ for
$0\!<\!t\!<\!1$ will then also be acceptable (see
Section~\ref{sec:assumptions}), but they will not satisfy the EE for
$\Btheta_0$, hence $\Bthetadot$ is not independent of the regressor
distribution for this conditional response distribution. ~~$\square$


\newpage

\subsection{Reweighting Diagnostics Applied to the Boston Housing Data}
\label{sec:reweighting-Boston}

\begin{figure}[h]
  \vspace{0in}
  \centerline{\includegraphics[width=5.35in]{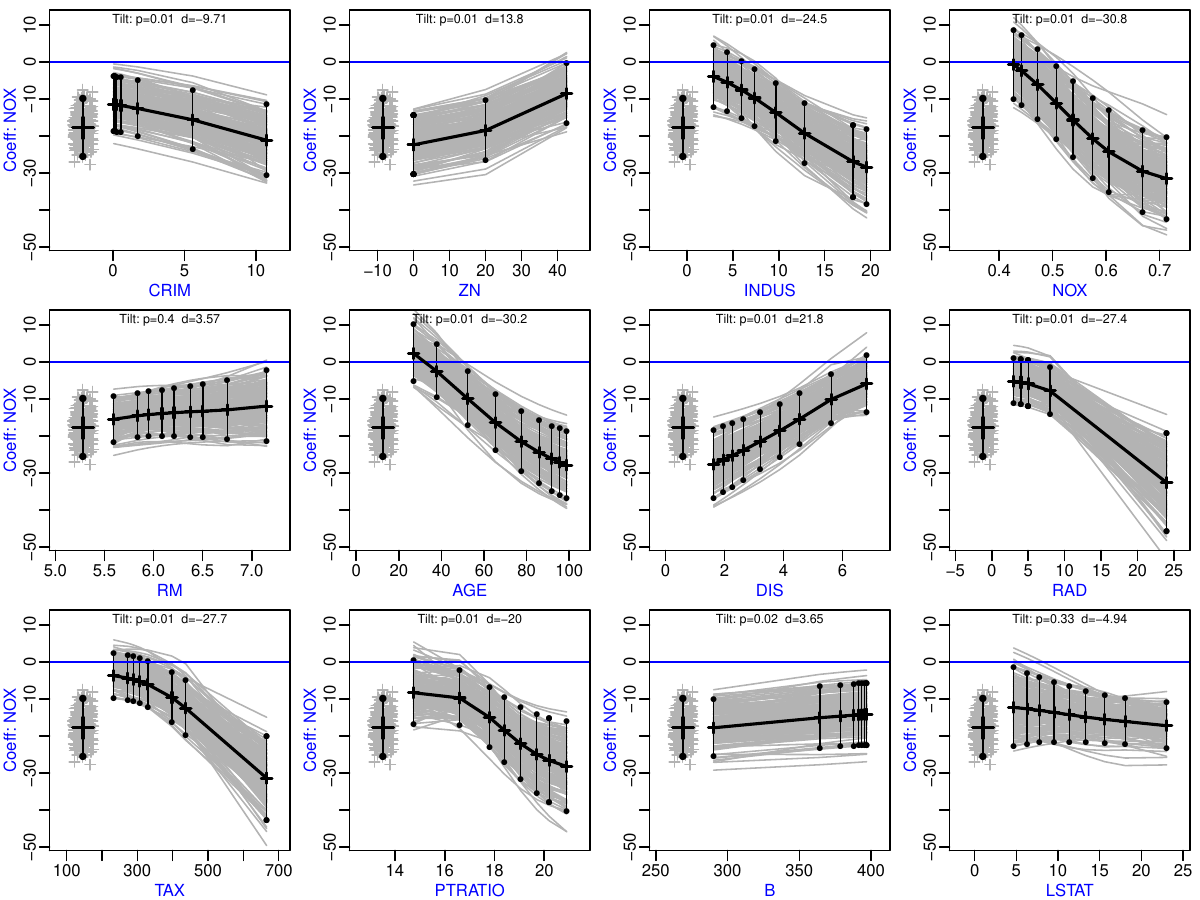}}
  \vspace{-.1in}
  \caption{\it Misspecification Diagnostics Based on Reweighting:
    Boston Housing Data.  The quantity of interest is the slope of the
    airpollution proxy {\tt NOX}.  In a hypothetical causal
    interpretation this slope is the ``effect size'' of {\tt NOX} on
    the response {\tt MEDV} (Median value of owner-occupied homes in
    \$1000's).  Misspecification of this quantity is evident, but in
    interpretable ways: The ``effect'' of {\tt NOX} is more negative
    on the values of homes, for example, at low values of {\tt DIS}
    (weighted distances to five Boston employment centres) and at high
    values of {\tt PTRATIO} (pupil-teacher ratio).  Also, the effect
    of {\tt NOX} is differentially greater at high values of {\tt
      NOX}, indicating an acceleration effect.  (The dummy regressor
    {\tt CHAS} for Charles River is left out because it is a binary
    variable with just 35 out of 506 cases at the high level and hence
    no non-trivial deciles for grid locations.)
  }\nocite{ref:BKY-2008}
  \label{fig:reweighting-Boston}
\end{figure}



\subsection{Reweighting and Partially Additive Models}
\label{sec:additive-models}

We discuss the type of diagnostic illustrated in
Section~\ref{sec:reweighting-on-own-regressor} where slopes
$\beta_j(\P)$ are reweighted wrt to their own regressors, $Z=X_j$.

If we interpret the slope $\beta_j(\P)$ of the regressor $X_j$ as a
partial derivative of the best approximation to the response surface,
then by localizing with a weight function $w_{\xi}(X_j)$ we
heuristically interpret $\beta_j(w_{\xi}(X_j)\P)$ as partial
derivative of the best approximation conditional on $X_j \approx \xi$.
If this localized partial derivative is not constant in $\xi$, it
indicates some form of nonlinearity of the response surface as a
function of $X_j$ (linearly adjusted for, and averaged over, all other
regressors).

This line of thinking may suggest that this diagnostic may be related
to a partially additive regression of the form
\begin{equation} \label{ref:partial-additive-model}
Y \approx s(X_j) \!+\! \sum_{k(\neq j)} \beta_k X_k ,
\end{equation}
where $s(\cdot)$ is a smooth function of $X_j$.  The heuristic
correspondence with the diagnostic is this:
\[
 s'(\xi) \approx \beta_j(w_{\xi}(X_j)\P).
\]
Note that in the diagnostic all adjustments are {\em linear} in the
other regressors, not nonlinearly additive, hence the diagnostic does
not correspond to a full additive model $Y \approx \sum_k s_k(X_k)$,
which adjusts $s_j(X_j)$ for all other $s_k(X_k)$.

Furthermore, there exists even a difference between the diagnostic and
the partially additive regression \eqref{ref:partial-additive-model}:
\begin{itemize} \itemsep 0em
\item The diagnostic adjusts $\beta_j(\cdot)$ linearly for all other $X_k$
  {\em with regard to the reweighted distribution}~$\, w_\xi(X_j) \P$, whereas
\item the partially additive regression
  \eqref{ref:partial-additive-model} adjusts $s(X_j)$ linearly for all
  other $X_k$ {\em with regard to the raw distribution} $\, \P$.
\end{itemize}
These subtle differences are to be considered when surprising effects
are observed under reweighting.

\subsection{A Connection of Reweighting to Nonparametrics}
\label{sec:connecting}

Generalizing the idea of reweighting and using many --- possibly
infinitely many --- weight functions provides a natural bridge from
parametric to nonparametric regression, namely, by using reweighting
functions that are ``running Parzen kernels'' as, for example, in
local linear smoothing.  The following are a few steps to describe the
general idea of localizing regression functionals when the regressor
space is $\Reals^p$: Let $\tw_\xvect(\xvec)$ be a family of Parzen
kernels, each member centered at a location $\xvect$ in regressor
space, an example being Gaussian kernels
$\tw_\xvect(\xvec) \sim \exp(-\|\xvec - \xvect\|^2/(2\sigma^2))$.
Then $w_\xvect(\xvec) = \tw_\xvect(\xvec)/\EP[w_\xvect(\Xvec)]$ is a
weight function that is normalized for $\PofX$ at each~$\xvect$.
Finally obtain the value of the regression functional localized
at~$\xvect$:
\begin{equation} \label{eq:reweighted-functionals}
  \Btheta_\xvect(\P) = \Btheta ( \, w_\xvect(\Xvec)\, \PofYandX \,) .
\end{equation}
Two special cases:
\begin{itemize} \itemsep 0.5em
\item If $\Btheta(\P) = \EP[Y]$, then
  $\xvect \mapsto \Btheta_\xvect(\P)$ is a regularized approximation
  to the response surface $\xvect \mapsto \EP[Y|\Xvec\!=\!\xvect]$,
  the result of local averaging.
\item If $\Bthetadot$ is the linear OLS functional, then
  $\Btheta_\xvect(\P)$ consists of a local intercept and local slopes
  at each location~$\xvect$, the latter forming a regularized
  approximation to the gradient at~$\xvect$.  If we define
  $f(\xvect) = \Btheta_\xvect(\P)' \xvect$, then $f(\xvect)$ is a
  locally linear approximation to the response surface
  $\xvect \mapsto \EP[Y|\Xvec\!=\!\xvect]$.
\end{itemize}
Estimating smooth functions and comparing them to linear ones has been
a diagnostic idea for some time, and a particularly useful approach
along these lines is by fitting additive models (Hastie and
Tibshirani, 1990)\nocite{ref:HT-1990}.  In the next subsection we will
pursue a different diagnostic idea that stays closer to the regression
functional of interest.


\begin{thebibliography}{9}





\bibitem{ref:Basu-et-al-1998} 
\textsc{Basu, A., Harris, I. R., Hjort, N. L.., Jones, M. C.} (1998).
Robust and Efficient Estimation by Minimising a Density Power Divergence.
\textit{Biometrika}
\textbf{85} (3), 549-559.


\bibitem{ref:Berk-et-al-2013}
\textsc{Berk, R., Brown, L., Buja, A., Zhang, K., and Zhao, L.} (2013).
Valid Post-Selection Inference.
\textit{The Annals of Statistics}
\textbf{41} (2), 802--837.



\bibitem{ref:BKY-2008}
\textsc{Berk, R. H.} and \textsc{Kriegler, B.} and \textsc{Yilvisaker, D.} (2008).
Counting the Homeless in Los Angeles County.
in \textit{Probability and Statistics: Essays in Honor of David A. Freedman},
Monograph Series for the Institute of Mathematical Statistics, D. Nolan and S. Speed (eds.)


\bibitem{ref:BGvZ-1997} %
\textsc{Bickel, P. J.} and \textsc{G\"otze, F.} and \textsc{van Zwet, W. R.} (1997).
Resampling Fewer than $n$ Observations: Gains, Losses, and Remedies for Losses.
\textit{Statistica Sinica}
\textbf{7}, 1--31.




\bibitem{ref:Breiman-1996} %
\textsc{Breiman, L.} (1996).
Bagging Predictors.
\textit{Machine Learning}
\textbf{24}, 123-140.









\bibitem{ref:BS-2016}
\textsc{Buja, A.} and \textsc{Stuetzle, W.} (2016, 2001).
Smoothing Effects of Bagging: Von Mises Expansions of Bagged Statistical Functionals.
\textit{arXiv:1612.02528}









\bibitem{ref:ET-1994} %
\textsc{Efron, B.} and \textsc{Tibshirani, R. J.} (1994).
\textit{An Introduction to the Bootstrap},
Boca Raton, FL: CRC Press.





\bibitem{ref:GP-2008}
\textsc{Gelman, A.} and \textsc{Park, D. K.} (2008).
Splitting a Regressor at the Upper Quarter or Third and the Lower Quarter or Third,
\textit{The American Statistician}
\textbf{62} (4), 1--8.

\bibitem{ref:GR-2007} 
\textsc{Gneiting, T.} and \textsc{Raftery, A. E.} (2007).
Strictly Proper Scoring Rules, Prediction and Estimation.
\textit{Journal of the American Statistical Association}
\textbf{102} (477), 359--378.


\bibitem{ref:Hall-1992} %
\textsc{Hall, P.} (1992).
\textit{The Bootstrap and Edgeworth Expansion.}
(Springer Series in Statistics)
New York, NY: Springer-Verlag.

\bibitem{ref:HRRS-1986} %
\textsc{Hampel, F. R.} and \textsc{Ronchetti, E. M.} and \textsc{Rousseeuw, P. J.} and \textsc{Stahel, W. A.} (1986).
\textit{Robust Statistics: The Approach based on Influence Functions.}
New York, NY: Wiley.

\bibitem{ref:Hansen-1982}
\textsc{Hansen, L. P.} (1982).
Large Sample Properties of Generalized Method of Moments Estimators.
\textit{Econometrica}
\textbf{50} (4), 1029-1054.

\bibitem{ref:HT-1990}
\textsc{Hastie, T. J.} and \textsc{Tibshirani, R. J.} (1990).
\textit{Generalized Additive Models},
London: Chapman \& Hall/CRC Monographs on Statistics \& Applied Probability.


\bibitem{ref:Boston-1978}
\textsc{Harrison, X.} and \textsc{Rubinfeld, X.} (1978).
Hedonic Prices and the Demand for Clean Air.
\textit{Journal of Environmental Economics and Management}
\textbf{5}, 81--102.


\bibitem{ref:Hausman-1978} 
\textsc{Hausman, J. A.} (1978).
Specification Tests in Econometrics.
\textit{Econometrica}
\textbf{46} (6), 1251-1271.



\bibitem{ref:Huber-1964} %
\textsc{Huber, P. J.} (1964).
Robust Estimation of a Location Parameter.
\textit{The Annals of Mathematical Statistics}
\textbf{35} (1) 73--101.

\bibitem{ref:Huber-1967} %
\textsc{Huber, P. J.} (1967).
The Behavior of Maximum Likelihood Estimation under Nonstandard Conditions.
\textsc{Proceedings of the Fifth Berkeley Symposium on Mathematical Statistics and Probability}, Vol. 1,
Berkeley: University of California Press, 221--233.

\bibitem{ref:IR-2015}
\textsc{Imbens, G. W.} and \textsc{Rubing, D. B.} (2015).
\textit{Causal Inference for Statistics, Social, and Biomedical Sciences: An Introduction.}
New York, NY: Cambridge University Press.



\bibitem{ref:Kent-1982} %
\textsc{Kent, J.} (1982).
Robust Properties of Likelihood Ratio Tests.
\textit{Biometrika}
\textbf{69} (1), 19--27.


\bibitem{ref:Kuchi-2018}
\textsc{Kuchibhotla, A. K} and \textsc{Brown, L. D.} and \textsc{Buja, A.} (2018).
Model-free Study of Ordinary Least Squares Linear Regression.
\textit{arXiv e-prints:1809.10538}.


\bibitem{ref:LZ-1986} %
\textsc{Liang, K.-Y.} and \textsc{Zeger, S. L.} (1986).
Longitudinal Data Analysis Using Generalized Linear Models.
\textit{Biometrika}
\textbf{73} (1), 13-22.









\bibitem{ref:Pearl-2009}
\textsc{Pearl, J.} (2009).
\textit{Causality: Models, Reasoning, and Inference (2nd ed.).}
New York, NY: Cambridge University Press.

\bibitem{ref:PBM-2016} %
\textsc{Peters, J.} and \textsc{B\"uhlmann, P.} and \textsc{Meinshausen, N.} (2016).
Causal Inference by Using Invariant Prediction: Identification and Confidence Intervals
\textit{Journal of the Royal Statistical Society, B}
\textbf{78}, Part~5, 947--1012.

\bibitem{ref:PR-1994} %
\textsc{Politis, D. N.} and \textsc{Romano, J. P.} (1994).
A General Theory for Large Sample Confidence Regions
based on Subsamples under Minimal Assumptions.
\textit{The Annals of Statistics}
\textbf{22}, 2031--2050.

\bibitem{ref:R-2008} %
\textsc{R Development Core Team} (2008).
R: A Language and Environment for Statistical Computing.
\textit{R Foundation for Statistical Computing,}
Vienna, Austria. ISBN 3-900051-07-0, URL http://www.R-project.org.


\bibitem{ref:Rieder-1994} %
\textsc{Rieder, H.} (1994).
\textit{Robust Asymptotic Statistics},
New York, NY: Springer-Verlag.



\bibitem{ref:Tukey-1962}
\textsc{Tukey,~J.~W.} (1962).
The Future of Data Analysis.
\textit{The Annals of Mathematical Statistics}
\textbf{33} (1), 1-67.




\bibitem{ref:White-1980a} %
\textsc{White, H.} (1980a).
Using Least Squares to Approximate Unknown Regression Functions.
\textit{International Economic Review}
\textbf{21} (1), 149-170.

\bibitem{ref:White-1980b} %
\textsc{White, H.} (1980b).
A Heteroskedasticity-Consistent Covariance Matrix Estimator and a Direct Test for Heteroskedasticity.
\textit{Econometrica}
\textbf{48}, 817-838.

\bibitem{ref:White-1981} %
\textsc{White, H.} (1981).
Consequences and Detection of Misspecified Nonlinear Regression Models.
\textit{Journal of the American Statistical Association}
\textbf{76} (374), 419-433.

\bibitem{ref:White-1982} %
\textsc{White, H.} (1982).
Maximum Likelihood Estimation of Misspecified Models.
\textit{Econometrica}
\textbf{50}, 1--25.

\bibitem{ref:White-1994} %
\textsc{White, H.} (1994).
\textit{Estimation, Inference and Specification Analysis.}
Econometric Society Monographs No.~22.
Cambridge, GB: Cambridge University Press.





\end{thebibliography}
\end{document}